\documentclass[12pt]{article}
\usepackage{times}                %   Postscript fonts
\usepackage{bm,amsbsy}
\usepackage{amsmath, amssymb,verbatim,graphicx, amsfonts}
\usepackage{mathrsfs}
\usepackage{verbatim}
\usepackage{multirow}
\usepackage{subfigure}
\usepackage{epsfig}
\usepackage{float}
\usepackage{pdfpages,setspace}
\usepackage[left=1in,top=1in,right=1in,bottom=1.0in]{geometry}
\usepackage{xcolor}
\usepackage{lineno}
\allowdisplaybreaks
\expandafter\def\expandafter\normalsize\expandafter{%
	\normalsize
	\setlength\abovedisplayskip{5pt}
	\setlength\belowdisplayskip{5pt}
	\setlength\abovedisplayshortskip{0pt}
	\setlength\belowdisplayshortskip{5pt}
}

\title{Free-stream preserving finite difference schemes for ideal magnetohydrodynamics on curvilinear meshes}

\author{
	Yize Yu
	\thanks{School of Mathematical Sciences, University of Science and Technology of China, Hefei, Anhui 230026, People's Republic of China
		{\tt yyz77@mail.ustc.edu.cn}. }
	\and
	Yan Jiang
	\thanks{School of Mathematical Sciences, University of Science and Technology of China, 
		Hefei, Anhui 230026, People's Republic of China
		{\tt jiangy@ustc.edu.cn}.  
		Research supported by NSFC grant 11901555.}
	\and	
	Mengping Zhang
	\thanks{School of Mathematical Sciences, University of Science and Technology of China, 
		Hefei, Anhui 230026, People's Republic of China
		{\tt mpzhang@ustc.edu.cn}. 
		Research supported by NSFC grant 11871448.}
}

\date{}

\begin{document}
	\maketitle

\begin{abstract}
	In this paper, a high order free-stream preserving finite difference weighted essentially non-oscillatory (WENO) scheme is developed for the ideal magnetohydrodynamic (MHD) equations on curvilinear meshes. 
	Under the constrained transport framework, magnetic potential evolved by a Hamilton-Jacobi (H-J) equation is introduced to control the divergence error. 
	In this work, we use the alternative formulation of WENO scheme \cite{christlieb2018high} for the nonlinear hyperbolic conservation law, and design a novel method to solve the magnetic potential. 
%	{\color{red}such that the scheme can preserve free-stream solutions of MHD equations}. 
%	In this work, we use the alternative formulation of WENO scheme \cite{christlieb2018high} for the nonlinear hyperbolic conservation law, and design a novel method which can maintain the linear function solution on both stationary and dynamically generalized coordinate systems for H-J equation to solve the magnetic potential.
%	In this work, a novel method is proposed to maintain the linear function solution on both stationary and dynamically generalized coordinate systems for H-J equation.
%	We use the alternative formulation of WENO scheme \cite{christlieb2018high} for the nonlinear hyperbolic conservation law, and the method proposed in this work to solve the magnetic potential. 
	Theoretical derivation and numerical results show that the scheme can preserve free-stream solutions of MHD equations, and reduce error
%	{\color{red}{
			more effectively than%}}
	 the standard finite difference WENO schemes for such problems.	
%	the proposed scheme can preserve free-stream solutions on both stationary and dynamically generalized coordinate systems, and reduce error effectively than the standard finite difference WENO schemes for such problems.
%	In this work, we use the alternative formulation of WENO scheme \cite{christlieb2018high} for the nonlinear hyperbolic conservation law, and design a novel method to solve the magnetic potential, such that the scheme can preserve free-stream solutions of MHD equations. 
%	In additional, metric terms are also carefully designed.
%	Theoretical derivation and numerical results show that the proposed scheme can preserve free-stream solutions on both stationary and dynamically generalized coordinate systems, and reduce error effectively than the standard finite difference WENO schemes for such problems.
\end{abstract}


	\section{Introduction}
	
	The ideal magnetohydrodynamic (MHD) equations are a fluid model to describe the dynamics of a perfectly conducting quasi-neutral plasma. 
	A range of areas including astrophysics and laboratory plasmas can be modeled with this system. 
	The magnetic field has to satisfy an extra constraint that its divergence is zero, which is a reflection of the principle that there are no magnetic monopoles.	
	Mathematically, the exact solution of the MHD equations preserves zero divergence for the magnetic field in future time, if the initial divergence is zero. 
%	{\color{red} 
		However, most standard numerical discretizations based on shock capturing methods do not propagate a discrete version of the divergence-free condition forward in time. The relevant information will be discussed later.%}
%	In this paper, we are interested in high order accurate conservative finite difference weighted essentially non-oscillatory (WENO) schemes for the ideal MHD system on curvilinear meshes (also referred to as generalized coordinate systems).}

	In this work, we are interesting in the free-stream preservation, which is an important property when finite difference schemes are applied to curvilinear meshes (also referred to as generalized coordinate systems). Failing preserving the free-stream condition may cause large errors and even lead to numerical instabilities for high-order schemes, see \cite{nonomura2010freestream, visbal2002use} for Euler equation and Navier-Stokes equation.
	This is because errors from nonpreserved free-stream may hide small physical oscillations, such as turbulent flow structures. 
%	Even though the standard finite difference weighted essentially non-oscillatory (WENO) schemes \cite{jiang1996efficient, balsara2000monotonicity, shu2009high} is widely employed for the simulation of complex flow fields due to their high order accuracy and good shock-capturing properties, and it has been extended to the ideal MHD equations successfully \cite{jiang1999high,  balsara2009divergence, shen2012cusp, balsara2013efficient, christlieb2014finite, christlieb2016high}.
%{\color{red}
	Here, we will consider the finite difference weighted essentially non-oscillatory (WENO) schemes. The standard WENO schemes \cite{jiang1996efficient, balsara2000monotonicity, shu2009high} are widely employed for the simulation of complex flow fields due to their high order accuracy and good shock-capturing properties. In particular, the methods have been studied to solve the ideal MHD equations successfully \cite{jiang1999high,  balsara2009divergence, shen2012cusp, balsara2013efficient, christlieb2014finite, christlieb2016high}.%}
However, it was found that standard finite difference WENO schemes could preserve the free-stream condition in the Cartesian coordinate system, but not in the generalized coordinate system 
%	 .In addition, the schemes in generalized grids have a rather large error arising from the metric terms on randomized grids and high dimensional wavy grids.This is 
due to non-exact cancellations of the metric derivatives when nonlinear reconstructions are performed. Details see \cite{cai2008performance, nonomura2010freestream}. 
More recently, Jiang et al. \cite{jiang2014free} studied an alternative formulation of finite-difference WENO scheme, originally proposed by Shu and Osher \cite{shu1988efficient}. 
Different from the standard one, in this formulation, the high order WENO interpolation procedure is applied to the solution, and the numerical technique for the metric terms \cite{visbal2002use} can be applied. 
	Moreover, theoretical derivation and numerical results showed that the finite difference WENO schemes based on the alternative flux formulation can preserve free-stream of Euler equation on both stationary and dynamically generalized coordinate systems. 
	Even though the alternative formulation could be used to solve ideal MHD equation on curvilinear meshes \cite{christlieb2018high}, however, the free-stream condition is not satisfied. The main reason is that the numerical method used to deal with zero divergence for the magnetic field may break the free-stream preserving property. 

	Controlling divergence errors in the magnetic field is one of the main numerical difficulties for simulating the ideal MHD equations. 
	Failure to control the divergence error in numerical computation could create an unphysical force parallel to the magnetic field (see \cite{brackbill1980effect} for instance), which may eventually result into numerical instabilities as its effects accumulate. 
	There are mainly four types of numerical approaches to solve this problem: the nonconservative eight-wave method \cite{gombosi1994axisymmetric},  the projection method \cite{brackbill1980effect}, the hyperbolic divergence cleaning method \cite{dedner2002hyperbolic}, and the various constrained transport methods \cite{evans1988simulation, dai1998divergence, balsara1999staggered,  fey2003constrained, balsara2004second,   rossmanith2006unstaggered, helzel2011unstaggered,  helzel2013high, christlieb2014finite, christlieb2018high}. See the review paper \cite{toth2000b} for more discussion.
	In \cite{christlieb2014finite, christlieb2018high}, authors used the unstaggered constrained transport method to address this issue. 
	In this framework, vector potential was introduced to correct the magnetic field. 
	Numerically, the conserved quantities of the ideal MHD equations and a magnetic vector potential are simultaneously evolved at the same time, with the evolution equations as hyperbolic conservation law and Hamilton-Jacobi (H-J) equation, respectively. 
	And then, set the magnetic field to be the (discrete) curl of the magnetic potential after each time step. 
	We find that for the free-stream flow, the magnetic potential should be a linear function in physical domain. 
	As a consequence, it is necessary to keep the linear solution exactly over time for the constrained transport equation on curvilinear coordinate systems numerically. 
	However, if we use the finite difference WENO method designed for H-J equations \cite{jiang2000weighted} on  lattice grid to solve magnetic potential on computational domain, it may fail to preserve the linear function. 
	Because the magnetic potential is not a linear function in computational domain any more.  One way to address this problem is that we can solve the equation on physical domain based on the WENO schemes for H-J equation on 2D triangular meshes \cite{zhang2003high}. 
	However, the algorithm employs the least square procedure and the computational cost would be very large for moving mesh problems.
	
	In this work, we will take the alternative formulation of finite difference WENO method to solve conservation laws, and develop a novel WENO method for the constrained transport equation or H-J equation, such that it can preserve the free-stream condition of ideal MHD equation on the generalized coordinate systems. 
	In the proposed scheme, we employ the standard WENO approximation to approach the derivatives on computational domain, coupled with the general Hamiltonian on unstructured meshes. Moreover, the metric terms are well designed to ensure they can be canceled exactly.
	
	The organization of the remaining sections is as follows. In Section 2, we will review the governing equations on generalized coordinate system. The numerical methods and the analysis of the free-stream-preservation property will be discussed in Section 3. In Section 4, extensive numerical examples are provided to demonstrate the free-stream preservation performance. Concluding remarks are given in Section 5.

\section{Governing equations}
	
In this section, we will briefly review the governing equations. 
Firstly, we start with the ideal MHD equations on Cartesian coordinates. 
A constrained transport framework is used here to control the divergence error of the magnetic field. 
And then, we will give the form of governing equations in curvilinear coordinates system.

\subsection{Ideal MHD equations}

The set of ideal MHD equations is a system of equations  describing the conservation of mass, momentum, energy, and magnetic field of a particular fluid. It is given in the following conservation form:
	\begin{equation}\label{mhd}
	\partial _t\left[
	\begin{array}{c}
	\rho \\
	\rho \textbf{u} \\
	E \\
	\textbf{B}
	\end{array}
	\right] + \nabla \cdot \left[
	\begin{array}{c}
	\rho \textbf{u}\\
	\rho \textbf{u} \otimes \textbf{u}+\left( p+\frac{1}{2}\|\textbf{B}\|^2\right) \mathbf{I}-\textbf{B} \otimes \textbf{B} \\
	\textbf{u} \left(E+p+\frac{1}{2}\|\textbf{B}\|^2\right) -\textbf{B}(\textbf{u} \cdot \textbf{B})\\
	\textbf{u} \otimes \textbf{B} - \textbf{B} \otimes \textbf{u}
	\end{array}
	\right]
	=0,
	\end{equation}
	
	\begin{equation}\label{div_free}
	\nabla \cdot \textbf{B}=0,
	\end{equation}	
	
	\noindent
	where, $\rho$ is the mass density, 
	$\rho \textbf{u}=\left(\rho u, \rho v, \rho w\right)^T$ is the momentum, 
	$E$ is the total energy,
	$\textbf{B}=\left(B_1, B_2, B_3\right)^T$ is the magnetic field,
	$p$ is the thermal pressure,
	and $\|\cdot\|$ is the Euclidean vector norm.
	The total energy density is given by
	\begin{equation*}\label{eq_of_state}
	E=\frac{p}{\gamma-1}
	+\frac{1}{2}\rho\|\textbf{u}\|^2
	+\frac{1}{2}\|\textbf{B}\|^2.
	\end{equation*}
	where $\gamma=5/3$ is the ideal gas constant.
	
	In particular, in the case of two dimensions considered in this work, all eight conserved variables, $\textbf{q}= ( \rho, \rho\textbf{u}, E, \textbf{B}) ^T$, can be non-zero, but each only depends on $x$, $y$, and $t$.  
		
	\subsection{Constrained transport}
	
	In order to reconstruct a discrete divergence-free magnetic field, we employ the constrained transport (CT) framework \cite{christlieb2014finite, christlieb2018high}. Instead of solving the magnetic field directly, a magnetic potential $\textbf{A}=(A_1, A_2, A_3)^T$ is introduced, satisfying 
	\begin{align}
	\label{eq:CT}
		\textbf{B}=\nabla \times \textbf{A}.
	\end{align}
	
	\noindent
	Furthermore, evolution equation for $\textbf{A}$ can be derived from the magnetic induction equation, as
\begin{equation}\label{ct}
\partial _t \textbf{A} + (\nabla \times \textbf{A})\times \textbf{u}=0.
\end{equation}	
See \cite{rossmanith2006unstaggered} for details.  
In the CT framework, we will simultaneously evolve the ideal MHD equations \eqref{mhd} and magnetic potential equations \eqref{ct} alongside each other, and set the magnetic field $\textbf{B}$ to be the (discrete) curl of the magnetic potential after each time step. 
Previous works \cite{christlieb2014finite, christlieb2018high} show such a procedure can control the divergence error in \textbf{B} and improve numerical stabilities of base schemes.

In the case of two dimensions, the divergence-free condition (\ref{div_free}) becomes	
\begin{equation*}
\nabla \cdot \textbf{B}=\partial _x B_1+\partial _y B_2=0.
\end{equation*}	
Thus $B_3$ has no impact on the divergence of $\textbf{B}$. It therefore suffices to only correct $B_1$ and $B_2$ in terms of controlling divergence errors. 
Furthermore, only the third component of $\textbf{A}$ need be evolved. 
Hence, for ease of presentation, we use a scalar quantity $A$ to denote the third component of $\textbf{A}$ in the following parts.
Consequently, the magnetic potential equations \eqref{ct} reduces to a scalar equation
\begin{equation}\label{ct1}
	\partial _t A +u\partial _x A+v\partial _y A=0	,
\end{equation}		
which is a H-J equation. While, $B_1$ and $B_2$ relate with $A$ by	
\begin{equation}
\label{eq:relation_AB}
	B_1=\partial _y A,\quad \text{and} \quad B_2=-\partial _x A.
\end{equation}

\subsection{Curvilinear coordinates and the free-stream condition}

\begin{figure}[htb]
	\centering
	\includegraphics[width=4.8in]{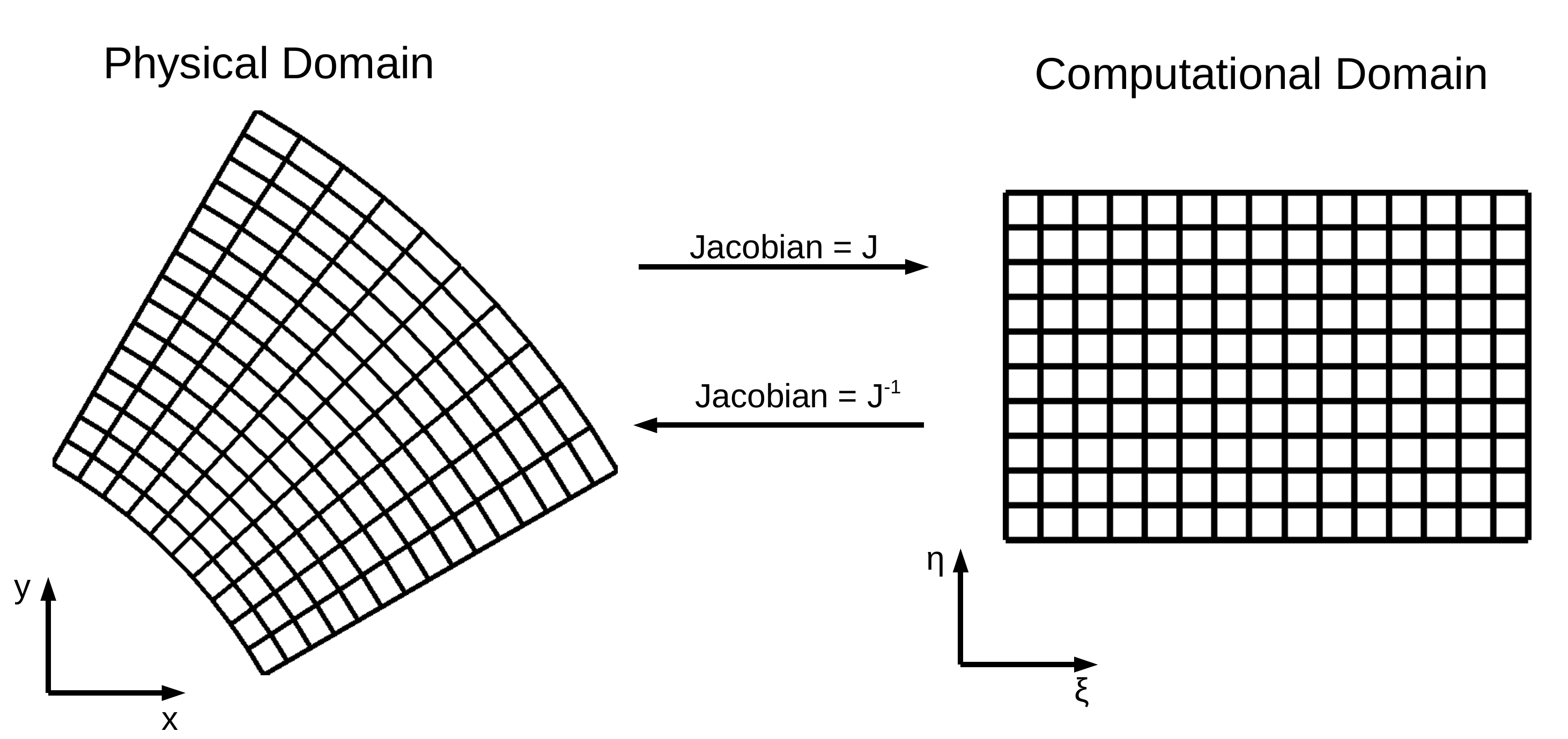}
	\caption{A schematic diagram of the transformations between the physical and computational domain.}
	\label{fig:CurvilinearMeshDomain}
\end{figure}

Here we provide a brief discussion of using curvilinear coordinates as	the computational domain to solve a general hyperbolic system and a H-J equation. 
Assume the Cartesian coordinates $\textbf{x}= (x,y,t)$ can be transformed in curvilinear coordinates $\textbf{r}=(\xi,\eta,\tau)$ via a continuous coordinate transformation. As illustrated in Fig. \ref{fig:CurvilinearMeshDomain}, a uniform mesh in the computational domain is typically used in our implementation.

Here, we choose $\tau=t$, and $J$ is the determinant of the Jacobian matrix defined by
\begin{equation}
\label{eq:Jacobian}
J^{-1}:=det \left[\frac{\partial \textbf{x}}{\partial \textbf{r}}\right]=
\left|\begin{array}{cccc} 
\partial _\xi x &    \partial _\eta x   & \partial _\tau x \\ 
\partial _\xi y &    \partial _\eta y   & \partial _\tau y \\ 
0 &    				0	 & 					1 
\end{array}\right|
=\partial _\xi x \,\partial _\eta y-\partial _\eta x \, \partial _\xi y.
\end{equation}
And the standard metrics satisfy
\begin{equation}\label{transform}
\begin{array}{lll}
	\partial _x \xi=J\, \partial _\eta y,
	&\partial _y \xi=-J\,\partial _\eta x,
	&\partial _t \xi=J \left( \partial _{\tau} y \, \partial _{\eta} x - \partial _{\tau} x \, \partial _{\eta} y \right), \\
	\partial _x \eta=- J\, \partial _\xi y,
	&\partial _y \eta=J\,\partial _\xi x, 
	&\partial _t \eta= J \left( \partial _{\tau} x \, \partial _{\xi} y -\partial _\tau y \, \partial _\xi x \right).	\\
\end{array}
%\left\{	
%\begin{aligned}	
%&\partial _x \xi=J\, \partial _\eta y,\\	
%&\partial _x \eta=- J\, \partial _\xi y,\\
%\end{aligned}
%\right.	
%\qquad 
%\left\{	
%\begin{aligned} 	
%&\partial _y \xi=-J\,\partial _\eta x,\\	
%&\partial _y \eta=J\,\partial _\xi x, \\
%\end{aligned}
%\right.	
%\qquad 
%\left\{	
%\begin{aligned}
%&\partial _t \xi=J \left( \partial _{\tau} y \, \partial _{\eta} x - \partial _{\tau} x \, \partial _{\eta} y \right), \\
%&\partial _t \eta= J \left( \partial _{\tau} x \, \partial _{\xi} y -\partial _\tau y \, \partial _\xi x \right).	
%%&\partial _t \xi=-\partial _\tau x \, \partial _x \xi-\partial _\tau y \, \partial _y \xi, \\
%%&\partial _t \eta=-\partial _\tau x \, \partial _x \eta-\partial _\tau y \, \partial _y \eta.
%\end{aligned}
%\right.
\end{equation}

For the system of hyperbolic conservation law in two dimensions,
\begin{equation*}\label{h-c-l-Cartesian}
\partial _t \textbf{q} + \partial _x \textbf{f}(\textbf{q}) + \partial _y \textbf{g}(\textbf{q})=0,
\end{equation*}
it can be transformed to the curvilinear coordinates, which has a conservative form 
\begin{equation}\label{h-c-l-curvilinear}
\partial _{\tau} \widetilde{\textbf{q}}+\partial _\xi \widetilde{\textbf{f}}+\partial _\eta \widetilde{\textbf{g}}=0.
\end{equation}
where,
\begin{equation*}
\widetilde{\textbf{q}}=J^{-1}\textbf{q}
,\qquad
\widetilde{\textbf{f}}=J^{-1}\left(\partial _t \xi \,\textbf{q}+\partial _x \xi \,\textbf{f} + \partial _y \xi \,\textbf{g} \right)
,\qquad
\widetilde{\textbf{g}}=J^{-1}\left(\partial _t \eta \,\textbf{q}+\partial _x \eta \,\textbf{f} + \partial _y \eta \,\textbf{g} \right).
\end{equation*}	

\noindent
When considering the magnetic potential equations \eqref{ct1} and \eqref{eq:relation_AB} in curvilinear mesh, they turn to 
\begin{align}
\label{eq:CT_curve}
	\partial _\tau A +  \left( u \partial _x \xi + v \partial _y \xi + J \partial_{t} \xi \right) \partial _\xi A + \left( u \partial _x \eta + v \partial _y \eta + J \partial_{t}\eta \right) \partial _\eta A  = 0,
\end{align}
and
\begin{align}
\label{eq:relation_AB_curve}
	B_1 =  \partial _y \xi \, \partial _\xi A + \partial _y \eta \, \partial _\eta A , \quad
	B_2 = -\partial _x \xi \, \partial _\xi A - \partial _x \eta \, \partial _\eta A ,
\end{align}
which are another H-J equation and the combination of partial derivatives of $A$, respectively.

In this work, we will look at the uniform flow (free-stream flow), i.e. $\textbf{q}=(\rho, \rho\textbf{u}, E, \textbf{B} )^T$ are constants. At this time, $\textbf{f}$ and $\textbf{g}$ are constants as well, and equation \eqref{h-c-l-curvilinear} can be simplified as
\begin{align*}
\partial_{\tau} \textbf{q} = -J \left( I_t \textbf{q} + I_x \textbf{f} + I_y \textbf{g}\right), 
\end{align*}
where,  
\begin{align}
\label{eq:SCL}
\left\{ \begin{array}{ll}
I_{t} =  \partial_{\tau} J^{-1} + \partial_{\xi} (J^{-1} \partial_{t} \xi) + \partial_{\eta} (J^{-1} \partial_{t}\tau),\\
I_{x} = \partial_{\xi} \left( J^{-1} \partial_{x}\xi \right) + \partial_{\eta} \left( J^{-1} \partial_{x}\eta \right), \\
I_{y} = \partial_{\xi} \left( J^{-1} \partial_{y}\xi \right) + \partial_{\eta} \left( J^{-1} \partial_{y}\eta \right). \\
\end{array}\right.
\end{align}
With the help of \eqref{eq:Jacobian} and \eqref{transform}, it is easy to check that $I_x = I_y = I_t =0$. 
Therefore, we have $\partial_{\tau} \textbf{q}=0$, that is, the uniform flow conditions are held. 
In a finite difference discretization, all three of the identities in equation \eqref{eq:SCL} must hold numerically to achieve the free-stream preserving condition. 
For stationary meshes, only the last two identities for $I_x$ and $I_y$ are required.

On the other hand, 
in the CT framework, $\mathbf{B}$ is determined through discreting $A$. 
We can prove that for the uniform flow, the initial magnetic potential should be in the form of $A(x,y,0) = B_1\, y - B_2\, x + C$, where $C$ is a constant.
In this case, the evolution equation \eqref{eq:CT_curve} tells us that 
\begin{align*}
	\partial_\tau A= (u - \partial_{\tau}x) B_2 - (v-\partial_{\tau}y) B_1.
\end{align*}
Therefore,
\begin{align*}
\label{eq:linear}
	A(\xi,\eta,\tau) = B_1\, y(\xi,\eta,\tau) - B_2\, x(\xi,\eta,\tau) + C + (u B_2 - v B_1) \tau .
\end{align*}
This means the magnetic potential $A$ is a linear function with respect to $(x,y)$ at any time.
Note that we set $\textbf{B}$ as the discrete curl of $A$.
Hence, in order to achieve the free-stream condition on curvilinear mesh numerically, we require the scheme can hold the linear solution of $A$ exactly.

Actually, it is easy to verify that the general H-J equation 
\begin{equation}\label{h-j}	
	\partial _t \phi +H(\partial _x \phi,\partial _y \phi) =0
\end{equation} 
would be transformed to the following one in the curvilinear coordinates 
\begin{equation}\label{h-j c}
	\partial _\tau \phi + \tilde{H}(\partial_{\xi}\phi, \partial_{\eta}\phi) =0,
\end{equation}	
where,
\begin{align}
	\tilde{H}(\partial_{\xi}\phi, \partial_{\eta}\phi)  = H(\{\partial _x \phi\},\{\partial _y \phi\})
	- \partial _\tau x \, \{\partial _x \phi\} - \partial _\tau y \, \{\partial _y \phi\},
\end{align}
and 
\begin{equation}	
\label{eq:curl_cureve}
	\{\partial _x \phi\}=
	\partial _x \xi \, \partial _\xi \phi + \partial _x \eta \, \partial _\eta \phi, 
	\qquad
	\{\partial _y \phi\}=
	\partial _y \xi \, \partial _\xi \phi + \partial _y \eta \, \partial _\eta \phi.
\end{equation}

\noindent
Moreover, it can preserve linear solution in physical domain, in the sense that if the initial condition is linear about $(x, y)$, i.e., $\phi(x,y,0)=C_1 x+C_2 y+C_3$ with constants $C_1, C_2$ and $C_3$, then the exact solution has the form that
$$ \phi(\xi,\eta,\tau)=C_1 x(\xi,\eta,\tau)+C_2y (\xi,\eta,\tau)-H(C_1,C_2) \tau+C_3. $$ 
Therefore, in the next section, we will design the high order finite difference WENO schemes, which can preserve the free-stream condition of conservation law \eqref{h-c-l-curvilinear}, or preserve linear function of H-J equation \eqref{h-j c} on curvilinear coordinate, respectively.

	\section{Numerical approach}
	
	In this section, we describe the construction of high order numerical methods to the governing equations.
	For time integration, we will use the third-order total-variation-diminishing (TVD) Runge-Kutta scheme for ordinary differential equation $u_\tau=L(u)$:
	\begin{equation}
	\label{eq:RK3}
	\begin{aligned}
	u^{(1)} =& u^n + \Delta \tau L(u^{n}),\\
	u^{(2)} =& \frac{1}{4} u^n + \frac{3}{4} \left( u^{(1)} +\Delta \tau L(u^{(1)}) \right),\\
	u^{n+1}=& \frac{1}{3} u^n + \frac{2}{3} \left( u^{(2)} +\Delta \tau L(u^{(2)}) \right).
	\end{aligned}
	\end{equation}

    \noindent
    Here, we will focus on the spatial discretization and the matrix terms. 
	Firstly, we give a brief description of the scheme for solving hyperbolic conservation systems, which can preserve the free-stream property.
	Then we deduce a novel numerical scheme for general H-J equations which can hold the linear solution exactly in arbitrary curvilinear meshes, even for moving mesh depending on time. 
	The computational domain is divided uniformly, with grid $\xi_i=i \Delta \xi$ and $\eta_j=j\Delta \eta$. 
		
\subsection{The numerical scheme for the ideal MHD equations}
	
	The ideal MHD equations in the curvilinear coordinates have a conservative form (\ref{h-c-l-curvilinear}), which is solved by a semi-discrete conservative finite difference scheme
	\begin{equation}\label{mhd_semi-discrete}
	\partial_{\tau} \widetilde{\textbf{q}}_{i,j}=
	- \frac{1}{\Delta \xi} \left( \widehat{\textbf{f}}_{i+1/2,j} - \widehat{\textbf{f}}_{i-1/2,j} \right) - 
	\frac{1}{\Delta \eta} \left( \widehat{\textbf{g}}_{i,j+1/2} - \widehat{\textbf{g}}_{i,j-1/2} \right).
	\end{equation}
	
	\noindent
	Here, $\widetilde{\textbf{q}}_{i,j}(\tau)$ is the numerical approximation to the point value $\widetilde{\textbf{q}}\left(\xi_{i}, \eta_{j}, \tau \right)$, and
	$\widehat{\textbf{f}}_{i+1/2,j}$ and $\widehat{\textbf{g}}_{i,j+1/2}$ are numerical fluxes sitting at the half points $\left(\xi_{i+1/2},\eta_{j}\right)$ and $\left(\xi_{i},\eta_{j+1/2}\right)$ respectively, satisfying
	\begin{subequations}
		\begin{align}
		\label{eq:order_conser_xi}
		& \frac{1}{\Delta \xi} \left( \widehat{\textbf{f}}_{i+1/2,j} - \widehat{\textbf{f}}_{i-1/2,j} \right) = \partial_{\xi} \widetilde{\textbf{f}}|_{i,j} + \mathcal{O}(\Delta \xi^k), \\
		\label{eq:order_conser_eta}
		& \frac{1}{\Delta \eta} \left( \widehat{\textbf{g}}_{i,j+1/2} - \widehat{\textbf{g}}_{i,j-1/2} \right)=\partial_{\eta} \widetilde{\textbf{g}}|_{i,j} + \mathcal{O}(\Delta \eta^k).
		\end{align}
	\end{subequations}
	
	In this work, we employ the numerical method introduced in \cite{christlieb2018high}, which based on an alternative formulation of the finite difference WENO scheme \cite{shu1988efficient, jiang2013alternative, jiang2014free}.
%	{\color{red}  
		In this framework, a high order WENO interpolation procedure is applied to the solution $\{\widetilde{\textbf{q}}\}$ rather than to the flux functions.%}
	For a multi-dimensional problem, the numerical fluxes can be obtained in a dimension-by-dimension fashion. Here, we take $\widehat{\textbf{f}}$ as an example,
	 \begin{equation}\label{flux}	
	 \widehat{\textbf{f}}_{i+1/2,j}=
	 \widetilde{\textbf{f}}_{i+1/2,j}
	 +\sigma_{i+1/2,j} \left(
	 -\frac{1}{24}\Delta \xi^2 \frac{\partial ^2\widetilde{\textbf{f}}}{\partial \xi^2}|_{i+1/2,j}+
	 \frac{7}{5760}\Delta \xi^4 \frac{\partial ^4\widetilde{\textbf{f}}}{\partial \xi^4}|_{i+1/2,j}\right),
	 \end{equation}
	 which gives fifth order accuracy in \eqref{eq:order_conser_xi}.
%	 {\color{red} 
	 	$\sigma_{i+1/2,j}$ is a filter to control numerical oscillations and will be defined later.%}
	 Likewise, the construction of numerical flux $\widetilde{\textbf{g}}_{i,j+1/2}$ can follow the similar idea in $\eta$-direction.

	 The first term of the numerical flux (\ref{flux}) is approximated by
	 \begin{equation}
	 \widetilde{\textbf{f}}_{i+1/2,j}=h \left( \textbf{q}^-_{i+1/2,j}, \textbf{q}^+_{i+1/2,j} \right),
	 \end{equation}
	 where, the values $\textbf{q}^{\pm}_{i+1/2,j}$ obtained by a WENO interpolation on $\textbf{q}$ in curvilinear coordinates $(\xi, \eta)$ in the $\xi$-direction for fixed $j$. Formulas of a fifth-order WENO interpolation are given in Appendix \ref{sec:WENO}.
	 The two-argument numerical function $h$ is an exact or approximate Riemann solver. In particular, the Lax-Friedrichs flux is used in this work
	 \begin{align*}
	 h \left( \textbf{q}^-, \textbf{q}^+ \right) = \frac{1}{2} \left( \, \widetilde{\textbf{f}}(\textbf{q}^-) + \widetilde{\textbf{f}}(\textbf{q}^+) -\alpha ( \textbf{q}^+  - \textbf{q}^-)\right), 
	 \end{align*}
	 where, $\alpha=\max_{1\leq l \leq 8} |\lambda_{l}(\textbf{q})|$ taken over the relevant range of $\textbf{q}$, and $\lambda_{l}(\textbf{q})$ are the eigenvalues of the Jacobian $\partial\widetilde{\textbf{f}} / \partial\textbf{q}$, which can be found in \cite{christlieb2018high}. 
	 The remaining higher order terms in (\ref{flux}) are constructed by central differences, 
	 \begin{equation}
	 \label{eq:high_order_terms}
	 \begin{aligned}
	 \Delta \xi^2 \frac{\partial ^2\widetilde{\textbf{f}}}{\partial \xi^2}|_{i+1/2,j}&=
	 \frac{1}{48} \left(
	 -  5 \, \widetilde{\textbf{f}}_{i-2,j} 
	 + 39 \, \widetilde{\textbf{f}}_{i-1,j}
	 - 34 \, \widetilde{\textbf{f}}_{i,j}
	 - 34 \, \widetilde{\textbf{f}}_{i+1,j}
	 + 39 \, \widetilde{\textbf{f}}_{i+2,j}
	 -  5 \, \widetilde{\textbf{f}}_{i+3,j}\right),\\
	 \Delta \xi^4 \frac{\partial ^4\widetilde{\textbf{f}}}{\partial \xi^4}|_{i+1/2,j}&=
	 \frac{1}{2} \left(
	 \widetilde{\textbf{f}}_{i-2,j} 
	 -3 \, \widetilde{\textbf{f}}_{i-1,j}
	 +2 \, \widetilde{\textbf{f}}_{i,j}
	 +2 \, \widetilde{\textbf{f}}_{i+1,j}
	 -3 \, \widetilde{\textbf{f}}_{i+2,j}
	 + \widetilde{\textbf{f}}_{i+3,j}\right).	 
	 \end{aligned}	 
	 \end{equation}
	 Both approximations give a truncation error of $\mathcal{O} (\Delta \xi^6)$, which guarantees a fifth-order accuracy of the numerical flux \eqref{flux}. 
	 
	Previous work \cite{jiang2013alternative, jiang2014free} showed that the numerical flux \eqref{flux} with $\sigma=1$ has a good performance of the compressible Euler equations in hydrodynamics. However, \cite{christlieb2018high} found that for the ideal MHD equations, an extra limiter $\sigma$ is necessary to switch the high-order numerical flux to a first-order flux near a strong discontinuity, since the linear approximations for high order terms \eqref{eq:high_order_terms} may cause oscillations . The filter $\sigma$ is carefully designed with the following property, 
	 \begin{align}
	 \sigma_{i+1/2,j}=
	 \left\{	\begin{array}{ll}
	 1+ O \left(\Delta \xi^3 \right), & \text{when \textbf{q} is smooth in the stencil $\{(\xi_{i-2},\eta_j), \ldots, (\xi_{i+3},\eta_j)\}$}, \\	
	 O \left(\Delta \xi^2 \right), & \text{when \textbf{q} contains a strong discontinuity in the stencil.}\\
	 \end{array} \right.	
	 \end{align}
	 Hence, the numerical flux \eqref{flux} can retain high-order accuracy in smooth regions, while control the oscillations around discontinuities.
	 The detailed formulation of $\sigma$ can be found in \cite{christlieb2018high}. 
	 
	 In particular, when $\textbf{q}$ are constants, we can easily prove that $\sigma=1$ always holds. In this case, the flux formula \eqref{flux} degenerates to the one used for the compressible Euler equations in hydrodynamics \cite{jiang2013alternative, jiang2014free}. Moreover, \cite{jiang2014free} showed that if we take the metrics at half points as follows,
	 \begin{equation} \label{eq:matric_linear}
	 \begin{aligned}
	 & \partial_{\xi}\gamma|_{i+1/2, j} = \frac{1}{256 \Delta \xi } \left( 3 \gamma_{i-2,j} -25\gamma_{i-1,j} +150\gamma_{i,j} + 150\gamma_{i+1,j} -25\gamma_{i+2,j} +3\gamma_{i+3,j}\right) , \\
	 & \partial_{\eta}\gamma|_{i, j+1/2} = \frac{1}{256 \Delta \eta} \left( 3 \gamma_{i,j-2} -25\gamma_{i,j-1} +150\gamma_{i,j} + 150\gamma_{i,j+1} -25\gamma_{i,j+2} +3\gamma_{i,j+3}\right) ,
	 \end{aligned}
	 \end{equation}
	 where $\gamma$ stands for $x$ and $y$, then finite difference WENO schemes with alternative formulation can preserve $I_x=I_y=0$ in \eqref{eq:SCL} exactly.  Hence, the proposed scheme \eqref{flux} can preserve free-stream scheme on stationary curvilinear mesh. 
	 
	 Furthermore, for moving meshes, we simply invoke the identity $I_t=0$ to evaluate $(J^{-1})_{\tau}$, i.e.,
	 \begin{align}
	 	(J^{-1})_{\tau}  =- \partial_{\xi} \left( \partial _{\tau} y \, \partial _{\eta} x - \partial _{\tau} x \, \partial _{\eta} y \right) 
	 	- \partial_{\eta} \left( \partial _{\tau} x \, \partial _{\xi} y -\partial _\tau y \, \partial _\xi x \right) .
	 \end{align}
	 Here, we use the linear scheme \eqref{eq:matric_linear} for spatial derivatives. And the temporal ones $(\partial _{\tau} x, \partial _{\tau} y)$ can be set as the exact derivatives or any numerical approach.   
	 Therefore, the free-stream conditions can be maintained in generalized curvilinear coordinate systems.
	 
\subsection{The numerical scheme for constrained transport equation}
	 In the previous work \cite{christlieb2018high}, the evolution equation \eqref{eq:CT_curve} is solved by a standard finite difference WENO method designed for H-J equations \cite{jiang2000weighted}, in which $A$ was treated as a function of $(\xi , \eta )$, and solved in a dimension-by-dimension fashion in the computational domain. Unfortunately, 
	 it can be proved that the standard WENO schemes could preserve the linear-function condition in the Cartesian coordinate system, but not in the generalized coordinate system. 
	 This {\color{red}is} because $A$ may be not linear for $\xi$ or $\eta$ in generalized coordinates, and the non-exact cancellations from nonlinear reconstructions will break the linear property.
	 
	 In the following part, we will propose a novel numerical schemes for H-J equation on curvilinear coordinates, in which we will approximate $\partial_x\phi$ and $\partial_y\phi$ in computational domain, and employ the Lax-Friedrichs type monotone Hamiltonian for arbitrary unstructured mesh developed by Abgrall in \cite{abgrall1996numerical}. In additional, the metrics are carefully designed such that the scheme can preserve the linear function exactly.

\subsubsection{Spatial discretization}	
\begin{figure}[htbp]        
		\center{\includegraphics[width=8cm] 
		%{fig9_new.png}} 
		{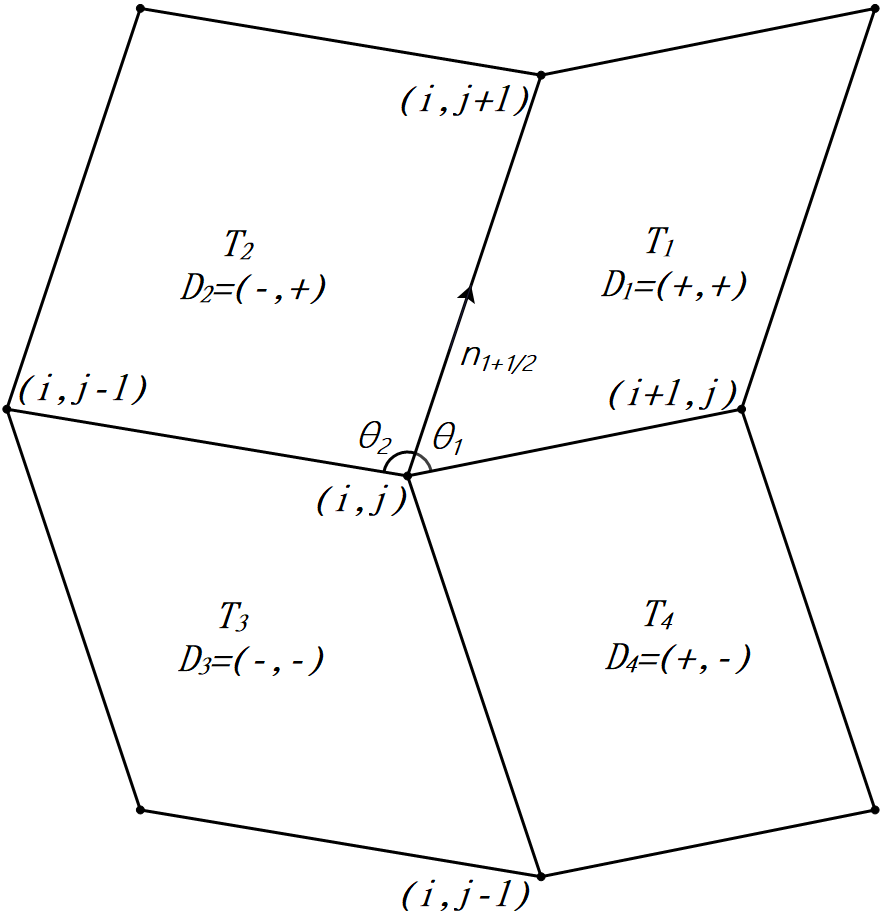}}       
		\caption{\label{Fig1} Node $(i,j)$ and its surrounding nodes. In this example, we have $\mathbf{D}_1=(+,+), \mathbf{D}_2=(-,+), 	\mathbf{D}_3=(-,-),  \mathbf{D}_4=(+,-)$.}      
	\end{figure}

	Let us look at the physical domain first. For simplify, we use node $(i,j)$ to stand for the point $\left( x_{i,j}, y_{i,j}\right)=\left( x\left(\xi_i,\eta_j\right), y\left(\xi_i,\eta_j\right)\right)$. 
	As shown in Fig. \ref{Fig1}, the area around node $(i,j)$ can be divided into four angular sectors, by the straight lines connecting $(i,j)$ and its surround nodes $(i\pm1,j)$ and $(i,j\pm1)$. 
	For each quadrilateral, denoted by $T_m$, $m=1,\ldots,4$, we use 
	%$D_m=(\pm,\pm)$
	$D_m=(\alpha,\beta)(\alpha=\pm,\beta=\pm)$  
	to indicate the positive/negative direction of $\xi$ and  $\eta$, or the quadrilateral contains the nodes $(i\pm 1, j)$ and $(i,j\pm 1)$. 
	For instance, $D_4=(+,-)$ means the region is enclosed with node $(i+1,j)$ and $(i,j-1)$. 
	In the follows, for each node $(i,j)$, we order the four quadrilaterals $T_{m}$ counterclockwise, starting from $D_{1}=(+,+)$. 
	See Fig. \ref{Fig1} as an example. 
	Additionally, we define $T_5=T_1$ for notational convenience. 
	
	We will denote by $\phi_{i,j}$ the numerical approximation to the solution at node $(i,j)$. 
	The equation \eqref{h-j c} can be solved by using a semi-discrete scheme
	\begin{align}
		\partial_\tau \phi_{i,j} + \widehat{H}\left( (\nabla \phi)^{1}_{i,j},(\nabla \phi)^{2}_{i,j},(\nabla \phi)^{3}_{i,j},(\nabla \phi)^{4}_{i,j} \right)=0.
	\end{align}
	Here, $(\nabla\phi)^{m}_{i,j}$ represents the numerical approximation of $\left(\{\partial _x\phi\},\{\partial _y\phi\}\right)$ (see \eqref{eq:curl_cureve}) at node $(i,j)$ in $T_m$.
	The monotone Hamiltonian for \eqref{h-j c} is given by
	\begin{equation}\label{monotone_Hamiltonian}
	\begin{aligned}
	\widehat{H}\left((\nabla \phi)^{1},(\nabla \phi)^{2},(\nabla \phi)^{3},(\nabla \phi)^{4}\right)
	=&H\left(\frac{1}{2\pi} \sum^{4}_{m=1}\theta_{m}(\nabla \phi)^{m} \right)
	-\left(x_\tau,y_\tau\right) \cdot \left(\frac{1}{2\pi} \sum^{4}_{m=1}\theta_{m}(\nabla \phi)^{m} \right)\\
	-&\frac{\lambda}{\pi}\sum^4_{m=1}\gamma_{m+\frac{1}{2}}
	\left(\frac{(\nabla \phi)^{m}+(\nabla \phi)^{m+1}}{2}\right)
	\cdot \textbf{n}_{m+\frac{1}{2}},
	\end{aligned}
	\end{equation}
	where,		
	$\theta_m$ is the inner angle of sector $T_m$,  
	and $\textbf{n}_{m+\frac{1}{2}}$ is the unit vector
	of the half-line $T_m \cap T_{m+1}$.
	The remaining variables are defined as follows,
	\begin{equation}\label{eq:parameter}
	\begin{aligned}	
	&\gamma_{m+\frac{1}{2}}=\tan(\frac{\theta_m}{2})+\tan(\frac{\theta_{m+1}}{2}),\\
	&\lambda = \max \left\{ \max_{\substack{A\leq u \leq B\\C\leq v \leq D}} \left|H_1(u,v)-x_\tau \right| , \max_{\substack{A\leq u \leq B\\C\leq v \leq D}} \left|H_2(u,v)-y_\tau \right|\right\}.
	\end{aligned}
	\end{equation}
	Here, $H_1$ and $H_2$ are the partial derivatives of $H$ with respect to $\partial _x \phi$ and $\partial _y \phi$, respectively. 
%	{\color{red} 
		And $[A, B]$ is the value range of $ \{\partial_x \phi\}^m_{i,j} $ and $[C, D]$ is the value range of $ \{\partial_y \phi\}^m_{i,j}$ for all $i,j,m$.%}

%	\begin{equation*}
%	\begin{aligned}
%	A=\min \left( \{\partial_x \phi\}^m_{i,j} \right), \quad  B=\max  \left( \{\partial_x \phi\}^m_{i,j} \right), \\
%	C=\min \left( \{\partial_y \phi\}^m_{i,j} \right), \quad D=\max \left( \{\partial_y \phi\}^m_{i,j} \right),
%	\end{aligned}
%	\end{equation*}

%	{\color{red}$A=\min \{\partial_x \phi\}^m_{i,j}$, $B=\max \{\partial_x \phi\}^m_{i,j}$ for all $i,j,m$.}
%	{\color{red}$C=\min \{\partial_y \phi\}^m_{i,j}$, $D=\max \{\partial_y \phi\}^m_{i,j}$ for all $i,j,m$.}
	%$[A,B]\times [C,D]$ is the value range for $(\nabla \phi)^{m}$.
	
	Next we will develop the high-order approximations for $(\nabla \phi)^{m}_{i,j}$. With the help of \eqref{eq:Jacobian} and \eqref{transform}, $\left(\{\partial _x\phi\},\{\partial _y\phi\}\right)$ can be rewritten as follows
	\begin{equation}
	\begin{aligned}
	\left(\{\partial _x\phi\},\{\partial _y\phi\}\right)
	=& \left( \partial _x \xi  \, \partial _\xi \phi  + \partial _x \eta \, \partial _\eta \phi, \, \partial _y \xi \, \partial _\xi \phi  + \partial _y \eta\, \partial _\eta \phi  \right) \\
	=& \left( \frac{\partial_\xi \phi \, \partial_\eta y -\partial_\eta \phi \, \partial_\xi y}
	{\partial_\xi x \, \partial_\eta y -\partial_\eta x \, \partial_\xi y}  , 
	 - \frac{\partial_\eta \phi \, \partial_\xi x -\partial_\xi \phi \, \partial_\eta x}
	{\partial_\xi x \, \partial_\eta y -\partial_\eta x \, \partial_\xi y}  \right).
	\end{aligned}
	\end{equation}
	Suppose $\mathbf{D}_m=(\alpha,\beta)$.  
	Then, we take
	\begin{align}
	(\partial _\xi \phi)^{m}_{i,j}= \partial _\xi \phi^{\alpha}_{i,j}, \quad 
	(\partial _\eta \phi)^{m}_{i,j}=\partial _\eta \phi^{\beta}_{i,j}.
	\end{align}
	Here, $\partial _\xi \phi^{\pm}_{i,j}$ and $\partial _\eta \phi^{\pm}_{i,j}$ are obtained by the 1D WENO approximation, based on the one-point left/right-biased stencils in $\xi$- and $\eta$-direction, respectively. Formulas with fifth order accuracy are given in Appendix \ref{append2}.

	In general, at node$(i,j)$, the schemes of $\partial_\xi \phi^{\pm}_{i,j}$ and $\partial_\eta  \phi^{\pm}_{i,j}$ have the following form,
	\begin{equation}
		\partial_\xi \phi^{\pm}_{i,j}= \sum^{q_1}_{l=-p_1} a^{\pm,l}_{i,j} \, \phi_{i+l,j},  \quad
		\partial_\eta \phi^{\pm}_{i,j}= \sum^{q_2}_{m=-p_2} b^{\pm,m}_{i,j} \, \phi_{i,j+m},
	\end{equation}
	with $\displaystyle \sum^{q_1}_{l=-p_1} a^{\pm,l}_{i,j}=0$ and $\displaystyle \sum^{q_2}_{m=-p_2} b^{\pm,m}_{i,j}=0$.
	If $\phi$ is linear for $x$ and $y$, with form
	\begin{equation}\label{linear_function}
		\phi=C_1 \, x+C_2 \, y-H(C_1,C_2) \, t+C_3,
	\end{equation}	
	and $C_1$, $C_2$ and $C_3$ are arbitrary constants, we can obtain that
	\begin{equation*}\label{scheme_for_linear}
	\begin{aligned}
		\partial_\xi \phi^\pm_{i,j}&= C_1 \, \sum^{q_1}_{l=-p_1} a^{\pm,l}_{i,j} \, x_{i+l,j}
		+C_2 \, \sum^{q_1}_{l=-p_1} a^{\pm,l}_{i,j} \, y_{i+l,j},\\	
		\partial_\eta \phi^\pm_{i,j}&= C_1 \, \sum^{q_2}_{m=-p_2} b^{\pm,m}_{i,j} \,x_{i,j+m}
		+C_2 \, \sum^{q_2}_{m=-p_2} b^{\pm,m}_{i,j} \,y_{i,j+m}.
	\end{aligned}
	\end{equation*}
	It's natural to notice that 
	$\displaystyle\sum^{q_1}_{l=-p_1} a^{\pm,l}_{i,j} \, k_{i+l,j}$ and $\displaystyle \sum^{q_2}_{m=-p_2} b^{\pm,m}_{i,j} \, k_{i,j+m}$ are high order approximations of $\partial_\xi k$  and $\partial_\eta k$ at node $(i,j)$ respectively, where $k$ stands for $x$ or $y$. Denote
	\begin{equation}\label{scheme_for_metrics}
	\begin{array}{lll}
	\displaystyle \widehat{\partial_\xi x}^{\pm}_{i,j}=\sum^{q_1}_{l=-p_1} a^{\pm,l}_{i,j} \, x_{i+l,j},&&
	\displaystyle \widehat{\partial_\xi y}^{\pm}_{i,j} =\sum^{q_1}_{l=-p_1} a^{\pm,l}_{i,j} \, y_{i+l,j},\\
	\displaystyle \widehat{\partial_\eta x}^{\pm}_{i,j}=\sum^{q_2}_{m=-p_2} b^{\pm,m}_{i,j} \, x_{i,j+m},&&
	\displaystyle \widehat{\partial_\eta y}^{\pm}_{i,j}=\sum^{q_2}_{m=-p_2} b^{\pm,m}_{i,j} \, y_{i,j+m}.
	\end{array}
	\end{equation}
	Therefore, we approximate the metric terms at $(i,j)$ in each $T_m$ by
	\begin{align}\label{eq:matric_form}
		(\partial_{\xi}k)^{m}_{i,j}=\widehat{\partial_\xi k}^{\alpha}_{i,j}, \quad 
		(\partial_{\eta}k)^{m}_{i,j}=\widehat{\partial_\eta k}^{\beta}_{i,j}.
	\end{align}
	This means on each grid point, we do the WENO method to obtain $\displaystyle\partial_\xi \phi^{\pm}_{i,j}$ and $\displaystyle\partial_\eta  \phi^{\pm}_{i,j}$ at first, and then use the nonlinear weights to simulate the metrics by  \eqref{scheme_for_metrics}. 

	Finally, we define the numerical approximations as 
	\begin{equation}\label{approximation2}
		\begin{aligned}
		\{\partial _x \phi\}^{m}_{i,j}
		&=\frac
		{\partial_\xi \phi^\alpha_{i,j} \, \widehat{\partial_\eta y}^\beta_{i,j} 
		-\partial_\eta \phi^\beta_{i,j} \, \widehat{\partial_\xi y}^\alpha_{i,j}}
		{\widehat{\partial_\xi x}^\alpha_{i,j} \, \widehat{\partial_\eta y}^\beta_{i,j}
		-\widehat{\partial_\eta x}^\beta_{i,j} \, \widehat{\partial_\xi y}^\alpha_{i,j}}, \quad
		\{\partial _y \phi\}^{m}_{i,j} = - \frac{\partial_\eta \phi^{\beta}_{i,j} \, \widehat{\partial_\xi x}^{\alpha}_{i,j} -\partial_\xi \phi^{\alpha}_{i,j} \, \widehat{\partial_\eta x}^{\beta}_{i,j} }
		{\widehat{\partial_\xi x}^\alpha_{i,j} \, \widehat{\partial_\eta y}^\beta_{i,j}
			-\widehat{\partial_\eta x}^\beta_{i,j} \, \widehat{\partial_\xi y}^\alpha_{i,j}} .
	\end{aligned}	
	\end{equation} 
	Moreover, when $\phi$ is given as \eqref{linear_function}, we can prove that
	\begin{equation*}
	\begin{aligned}
		\{\partial _x \phi \}^{m}_{i,j}
		&=\frac
		{\left(C_1 \, \widehat{\partial_\xi x}^{\alpha}_{i,j} +C_2 \, \widehat{\partial_\xi y}^{\alpha}_{i,j}\right)  \, \widehat{\partial_\eta y}^{\beta}_{i,j} - \left(C_1 \, \widehat{\partial_\eta x}^{\beta}_{i,j}+C_2 \, \widehat{\partial_\eta y}^{\beta}_{i,j}\right) \, \widehat{\partial_\xi y}^{\alpha}_{i,j}}
		{\widehat{\partial_\xi x}^\alpha_{i,j} \, \widehat{\partial_\eta y}^\beta_{i,j}
		-\widehat{\partial_\eta x}^\beta_{i,j} \, \widehat{\partial_\xi y}^\alpha_{i,j}}
		=C_1. \\
		\{\partial _y \phi\}^{m}_{i,j} 
		&= - \frac{ \left(C_1 \, \widehat{\partial_\eta x}^{\beta}_{i,j}+C_2 \, \widehat{\partial_\eta y}^{\beta}_{i,j}\right)  \, \widehat{\partial_\xi x}^{\alpha}_{i,j} -\left(C_1 \, \widehat{\partial_\xi x}^{\alpha}_{i,j} +C_2 \, \widehat{\partial_\xi y}^{\alpha}_{i,j}\right)  \, \widehat{\partial_\eta x}^{\beta}_{i,j} }
		{\widehat{\partial_\xi x}^\alpha_{i,j} \, \widehat{\partial_\eta y}^\beta_{i,j}
			-\widehat{\partial_\eta x}^\beta_{i,j} \, \widehat{\partial_\xi y}^\alpha_{i,j}}  
		= C_2.
	\end{aligned}	
	\end{equation*} 
	Therefore, according to the consistence of this monotone Hamiltonian \cite{abgrall1996numerical}, we have
	\begin{equation}\label{linear_condtion_spatial}
	\widehat{H}\left((\nabla \phi)^{1},(\nabla \phi)^{2},(\nabla \phi)^{3},(\nabla \phi)^{4}\right)
	=H(C_1,C_2)-(C_1,C_2) \cdot(\partial _\tau x_{i,j} ,\partial _\tau y_{i,j}).
	\end{equation}
	It's obviously that condition (\ref{linear_condtion_spatial}) is enough to preserve linear solution exactly for stationary meshes ($\partial_{\tau}x=\partial_{\tau}y=0$). 
	But for moving meshes, the discrete schemes of $\partial _\tau x$ and $\partial _\tau y$ should be consistent with time integration scheme. 
	
	\subsubsection{Temporal derivatives}
	For simplicity, we will start with the Euler forward method at first, 
	\begin{align}
	\phi_{i,j}^{n+1} = \phi_{i,j}^{n} - \Delta \tau \, \widehat{H}\left( (\nabla \phi)^{1}_{i,j}, (\nabla \phi)^{2}_{i,j}, (\nabla \phi)^{3}_{i,j}, (\nabla \phi)^{4}_{i,j} \right).
	\end{align}
	
	\noindent 
	In particular, suppose $\phi$ satisfies the linear solution (\ref{linear_function}) at $\tau=\tau_n$ exactly,
	$$\phi_{i,j}^{n}=\phi(x^{n}_{i,j},y^{n}_{i,j},\tau_n)=-\tau_n H(C_1,C_2)+(C_1,C_2) \cdot (x^{n}_{i,j},y^{n}_{i,j}) + C_3.$$
	To ensure this is still true for $\phi^{n+1}_{i,j}$
	\begin{align*}
	\phi_{i,j}^{n+1}
	&=\phi(x^{n+1}_{i,j},y^{n+1}_{i,j},\tau_{n+1})=-\tau_{n+1} H(C_1,C_2)+(C_1,C_2) \cdot (x^{n+1}_{i,j},y^{n+1}_{i,j}) + C_3\\
	&=-\tau_n H(C_1,C_2)+(C_1,C_2) \cdot (x^{n}_{i,j},y^{n}_{i,j}) + C_3 \\
	& + \Delta \tau\left(- H(C_1,C_2) + (C_1,C_2) \cdot (\frac{x^{n+1}_{i,j}-x^{n}_{i,j}}{\Delta \tau}, \frac{y^{n+1}_{i,j}-y^{n}_{i,j}}{\Delta \tau}) \right),
	\end{align*}	
	
	\noindent
	we can set the time derivative $(x_\tau,y_\tau)$ at $\tau_n$ as
	\begin{equation}
	\begin{aligned}	
	(\widehat{\partial _\tau x}^n_{i,j} ,\widehat{\partial _\tau y}^n_{i,j}) = ( \frac{x(\xi_{i},\eta_{j},\tau_{n+1})-x(\xi_{i},\eta_{j},\tau_{n})}{\Delta\tau}, \frac{y(\xi_{i},\eta_{j},\tau_{n+1})-y(\xi_{i},\eta_{j},\tau_{n})}{\Delta\tau} ).	
	\end{aligned}	 	
	\end{equation}
	Extension to the third-order TVD Runge-Kutta method \eqref{eq:RK3} is straightforward, since the schemes are convex combinations of Euler forward steps.
	
	\begin{equation}
	\begin{aligned}
	(\widehat{\partial _\tau x}^n_{i,j} ,\widehat{\partial _\tau y}^n_{i,j}) =& \frac{1}{\Delta \tau}\left(x_{i, j}^{(1)}-x_{i, j}^{n}  ,
	y_{i, j}^{(1)}-y_{i, j}^{n}\right), \\ 
	(\widehat{\partial _\tau x}^{(1)}_{i,j} ,\widehat{\partial _\tau y}^{(1)}_{i,j}) =& \frac{1}{\Delta \tau}
	\left( 4 x_{i, j}^{(2)}-3 x_{i, j}^{n} - x_{i, j}^{(1)} ,  4 y_{i, j}^{(2)}-3 y_{i, j}^{n} - y_{i, j}^{(1)} \right), \\
	(\widehat{\partial _\tau x}^{(2)}_{i,j} ,\widehat{\partial _\tau y}^{(2)}_{i,j}) =& \frac{1}{\Delta \tau}
	\left( 3 x_{i, j}^{n+1} - x_{i, j}^{n} - 2 x_{i, j}^{(2)} , 3 y_{i, j}^{n+1} - y_{i, j}^{n} - 2 y_{i, j}^{(2)} \right),
	\end{aligned}
	\end{equation}
	where, $\gamma_{i, j}^{(1)}=\gamma(\xi_{i},\eta_{j},\tau_{n}+\Delta \tau)$, $\gamma_{i, j}^{(2)}=\gamma(\xi_{i},\eta_{j},\tau_{n}+1/2\Delta \tau)$, and $\gamma$ stands for $x$ or $y$.

\subsection{Algorithm}
	We summarize this framework below.
	We write the semi-discrete form of MHD equations (\ref{mhd_semi-discrete}) as
	\begin{equation}\label{semi_dis_mhd}
	\widetilde{\textbf{q}}'\left(\tau\right)
	=\mathcal{L} \left( \widetilde{\textbf{q}} \left(\tau\right) \right),
	\end{equation}
	where $\widetilde{\textbf{q}}=J^{-1} \left(\rho,\rho\textbf{u},E,\textbf{B}\right)$.
	And the evolution equation of the magnetic potential (\ref{ct}) in semi-discrete form is 
	\begin{equation}\label{semi_dis_ct}
	A\left(\tau\right)
	=\mathcal{H}\left(A\left(\tau\right),\textbf{u}\left(\tau\right)\right).
	\end{equation}
	A single time-step from time $\tau=\tau_n$ to time $\tau=\tau_{n+1}$ with Euler forward consists of the following sub-steps:
	
	\begin{itemize}
		\item [(0)] 
		Start with $\widetilde{\textbf{q}}^n$ and $A^{n}$ (the solutions at $\tau _n$). 
		     
		\item [(1)]
		Build the right-hand sides of systems (\ref{semi_dis_mhd})  and (\ref{semi_dis_ct}) with WENO methods. 
%		{\color{red} 
			Note that the metric terms in $\mathcal{L}\left(\widetilde{\textbf{q}}^n\right)$ are obtained by \eqref{eq:matric_linear}, while the terms in $\mathcal{H}\left(A^n,\textbf{u}^n\right)$ are obtained by \eqref{scheme_for_metrics}. %}
		
		\item[(2)] Update each system independently:
		\begin{align*}
		\widetilde{\textbf{q}}^{*}
		&=\widetilde{\textbf{q}}^n+\Delta \tau \mathcal{L}\left(\widetilde{\textbf{q}}^n\right),\\
		A^{n+1}
		&=A^{n} +\Delta \tau \mathcal{H}\left(A^n,\textbf{u}^n\right),	 		
		\end{align*}
		where $\textbf{q}^*=1/J^{n+1} \left(\rho^{n+1},\rho\textbf{u}^{n+1},E^*,\textbf{B}^*\right)$.
		$\textbf{B}^*$ is the predicted magnetic field that in general does not satisfy a discrete divergence-free constraint and $E^*$ is the predicted energy.
		
		\item [(3)]
		Replace $\textbf{B}^*$ by a discrete curl of the magnetic potential $A^{n+1}$:
		\begin{equation}\label{eq:discrete}
		\begin{aligned}
		(B_1)^{n+1}_{i,j} =&  J^{n+1}_{i,j} \left( -( \partial _\eta x )^{n+1}_{i,j} \, (\partial _\xi A)^{n+1}_{i,j}  + ( \partial _\xi x)^{n+1}_{i,j}  \, (\partial _\eta A)^{n+1}_{i,j} \right) , \\
		(B_2)^{n+1}_{i,j}  =& J^{n+1}_{i,j} \left( -( \partial _\eta y)^{n+1}_{i,j}  \, (\partial _\xi A)^{n+1}_{i,j}  + ( \partial _\xi y)^{n+1}_{i,j}  \, (\partial _\eta A)^{n+1}_{i,j} \right) .
		\end{aligned}
		\end{equation}
		Here, we only use the sixth-order central differences in the constrained transport step
		\begin{equation}
	\begin{aligned}
	\partial_{\xi} k_{i,j} \approx \frac{1}{60\Delta \xi}	\left(-  k_{i-3, j}+9  k_{i-2, j} -45 k_{i-1, j} +45 k_{i+1, j} -9 k_{i+2, j} +  k_{i+3, j} \right), \\ 
	\partial_{\eta} k_{i,j}\approx  \frac{1}{60 \Delta \eta} \left(- k_{i, j-3} +9 k_{i, j-2} -45 k_{i, j-1} +45 k_{i, j+1} -9 k_{i, j+2} + k_{i, j+3}\right),
	\end{aligned}
	\end{equation}
	where $k$ stands for $A$, $x$ and $y$.
		
		\item [(4)]	
		Set the corrected total energy density $E^{n+1}$. Here, we keep the pressure the same before and after the magnetic field correction step:
		\begin{align*}
		E^{n+1}= E^{*}+\frac{1}{2}\left(\Vert{\textbf{B}^{n+1}}\Vert ^2-\Vert\textbf{B}^*\Vert ^2 \right),
		\end{align*}
		with $\|\textbf{B}\|^2 = B_1^2 + B_2^2 + B_3^2$.
	\end{itemize}

We want to remark that when $\phi$ is the linear function \eqref{linear_function}, the discretization \eqref{eq:discrete} can give $B_1$ and $B_2$ exactly. On the other hand, we note that while such discretization only guarantees the divergence free condition of magnetic field to truncation errors on general curvilinear meshes. In \cite{christlieb2018high}, authors found this is sufficient to suppress the unphysical oscillations associated with the divergence error of $\textbf{B}$.

	\section{Numerical results}
	In this section, the numerical simulations for both Hamilton-Jacobi  equations and the ideal MHD equations will be presented.
	
	At first, we will study the performance of proposed scheme in solving H-J equations on wavy mesh, comparing with the classical WENO scheme \cite{jiang2000weighted}. 
	For simplicity, the scheme we proposed in this paper is denoted as PL-WENO, and the classical WENO scheme is denoted as NPL-WENO. 
	Next, we will discuss the results of ideal MHD equations coupled with constrained transport equation on both stationary and dynamical meshes.  
	For comparison, we also test the method in \cite{christlieb2018high}. 
	The only difference between those two methods is the scheme solving the constrained transport equation, that our method uses PL-WENO, while \cite{christlieb2018high} used the classical NPL-WENO, which can not preserve free-stream condition.
	
    Here, We only test with the fifth order WENO schemes in space. 
    And the third order TVD Runge-Kutta scheme \eqref{eq:RK3} is used for time integration.
    The computational domain is taken as $(\xi,\eta)\in[0,L_x]\times[0,L_y]$ without special declaration, and take the mesh sizes
    $$\Delta \xi =L_x/(I_{max}-1) ,\quad \Delta\eta = L_y/(J_{max}-1).$$
	
	\subsection{H-J Equation: accuracy test} \label{exadd1}
	We consider the problem that 
	\begin{equation*}
	\left\{
	\begin{array}{lr}
	\phi_{t}=\phi_x+\phi_y, \\
	\phi(x,y,0)=\sin(x+y). \\		
	\end{array}
	\right.
	\end{equation*}
	on a stationary wavy grid
	\begin{equation}\label{wavy_grid}
	\begin{aligned}
	x_{i,j}&=x_{min}+\Delta \xi \left[ (i-1)+A_x \sin{\frac{N_y(j-1)\Delta \eta}{L_y}} \right],\\
	y_{i,j}&=y_{min}+\Delta \eta \left[ (j-1)+A_y \sin{\frac{N_x(i-1)\Delta \xi}{L_x}} \right],		
	\end{aligned}
	\end{equation}
	where,
	\begin{equation*}
	\begin{aligned}
	&i=1,2,\cdots,I_{max},\qquad
	j=1,2,\cdots,J_{max},\\
%	&\Delta x_{0}=\frac{L_x}{I_{max}-1},\qquad 
%	\Delta y_{0}=\frac{L_y}{J_{max}-1},\\	
	&x_{min}=-\frac{L_{x}}{2},\qquad
	y_{min}=-\frac{L_{y}}{2}.
	\end{aligned}
	\end{equation*}
	The wavy grid parameters used in this test are $L_x=L_y=2\pi$, $\Delta\xi A_x=0.01$, $\Delta \eta A_y=-0.02$, $N_x=N_y=2\pi$. And we test with $I_{max} = J_{max} = 41, 81, 161$ and $321$ successively.
	To realize fifth order in time, the time step is choose as $\Delta \tau \approx \Delta \xi^{5 / 3}$.
	The boundary condition are all periodic	for this smooth test.
	The $L_1$ and $L_\infty$ errors and orders of accuracy at time $T=0.5$ are presented in Table \ref{accuracy}, showing that the scheme can achieve the designed fifth order accuracy.	
	
	\begin{table}
		\caption{Example \ref{exadd1}: $L_1$ and $L_\infty$ errors at $T=0.5$}
		\bigskip
		\centering
		\begin{tabular}{ccccc}
			\hline
			$I_{max} \times J_{max}$ & $L_1$ errors & order & $L_\infty$ errors & order \\  \hline
			$41 \times 41$ & 9.37E-06 &  --   & 2.05E-05 &  --  \\
			$81 \times 81$ & 2.96E-07 & 4.89 & 6.08E-07 & 5.08\\
			$161\times161$ & 9.21E-09 & 5.00 & 1.84E-08 & 5.04\\
			$321\times321$ & 2.86E-10 & 5.01 & 5.67E-10 & 5.02\\
			\hline
		\end{tabular}
		\label{accuracy}
	\end{table}

	\subsection{H-J Equation: 2D nonconvex Riemann problem}\label{ex0}
\begin{figure}[htbp]
		\centering	
		 \subfigure[Grid]{
		\begin{minipage}[t]{0.33\linewidth}
			\centering
			\includegraphics[width=2in]{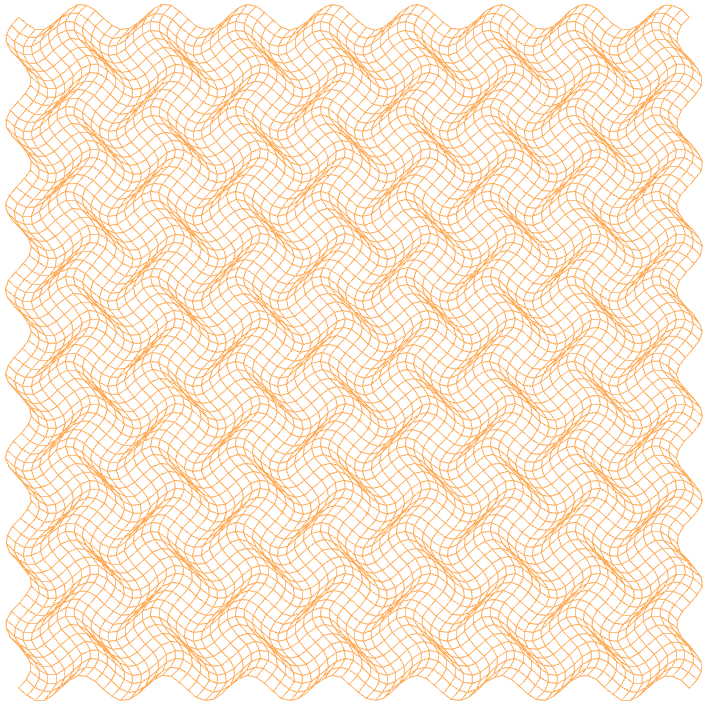}
		\end{minipage}%
		}%
	   \subfigure[PL-WENO]{
			\begin{minipage}[t]{0.33\linewidth}
				\centering
				\includegraphics[width=2in]{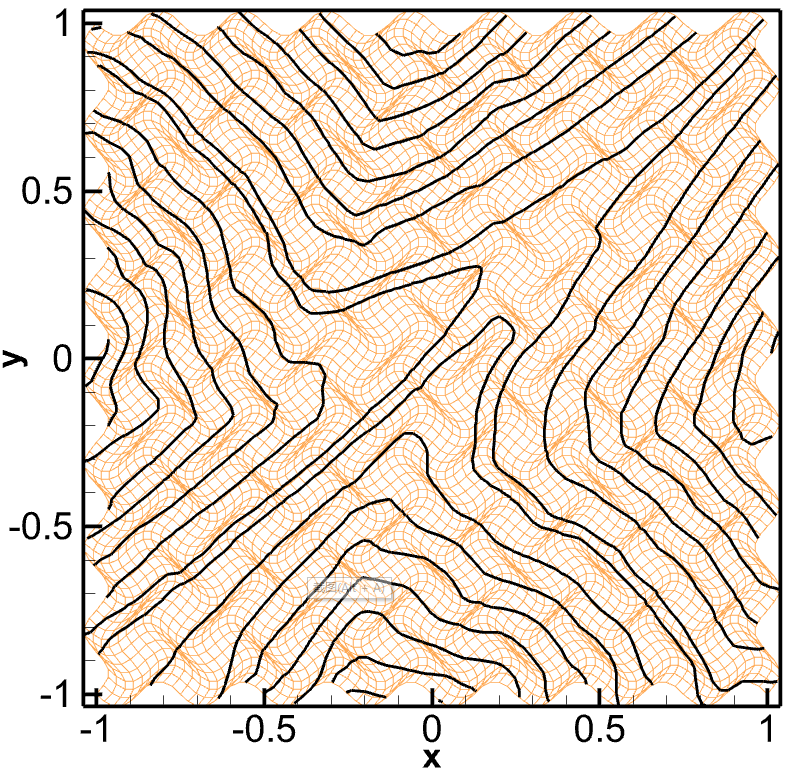}
			\end{minipage}%
		}%
		\subfigure[NPL-WENO]{
			\begin{minipage}[t]{0.33\linewidth}
				\centering
				\includegraphics[width=2in]{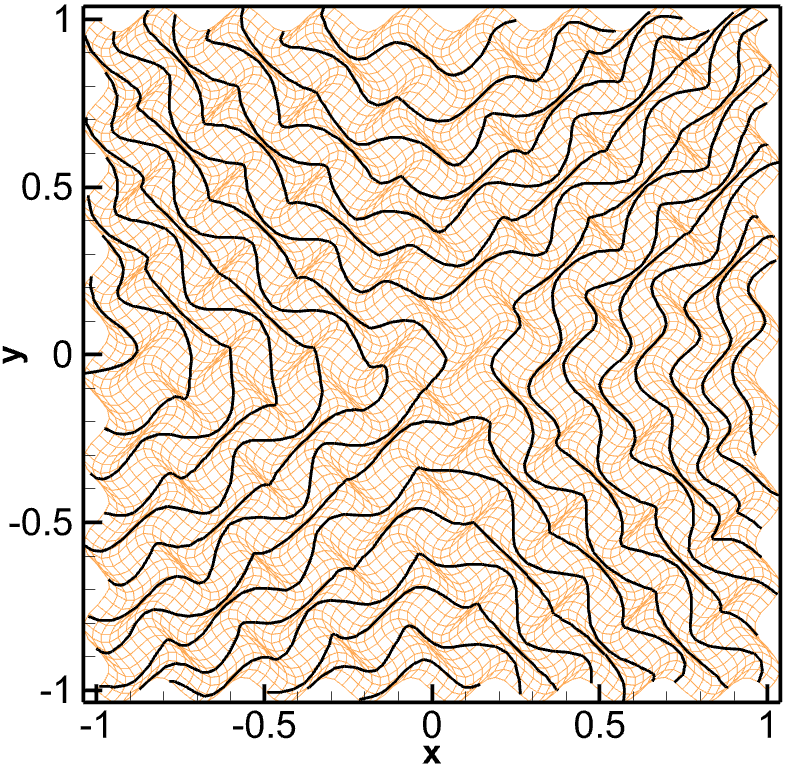}
			\end{minipage}%
		}%		  
		\centering
		\caption{Example \ref{ex0}: 16 contour lines of $\phi$ between $-2.5$ and $2.5$ are shown.}
		\label{example0}
	\end{figure}

	We solve a 2D nonconvex Riemann problem from \cite{osher1991high},
	\begin{equation*}
	\left\{
	\begin{aligned}
	&\partial_t \phi + \sin\left(\partial_x \phi+\partial_y \phi \right)=0,\\
	&\phi (x,y,0)=\pi(\left|y\right|-\left|x\right|),
	\end{aligned}
	\right.	
	\end{equation*}
	on a stationary wavy grid expressed as \eqref{wavy_grid}.
	The wavy grid parameters used in this test are $I_{max}=J_{max}=81$, $L_x=L_y=2$, $A_x=1.5$, $A_y=-1.5$, $N_x=N_y=16\pi$.
	The time step is chosen as $\Delta \tau= CFL \,\Delta \xi / \lambda$,
	with $CFL=0.1$, and $\lambda$ is given in (\ref{eq:parameter}).
	We compute the solutions up to time $T=1$. Numerical solutions are plotted in Fig. \ref{example0}.
	It is observed that the PL-WENO can keep the solution's shape from the influence of the wavy grid and give better simulation to the viscosity solutions. \cite{osher1991high}, while NPL-WENO bend the linear part of solution we do not want.
	
	\subsection{MHD: Free-stream preserving properties}\label{ex1}
	In this example, an free-stream fluid is imposed, with $\rho = \gamma^2$, $p= \gamma$, $\textbf{u}=(1,0,0)^T$, and $\textbf{B}=(1,1,0)^T$. Thus the $y$-direction velocity $v$ and the $z$-direction velocity $w$ are expected to remain machine zero. 
	Here, we test with four different grids:
	\begin{itemize}
		\item \textbf{Stationary wavy grid}. It is expressed as \eqref{wavy_grid} with wavy grid parameters $I_{max}=J_{max}=41$, $L_x=L_y=4\pi$, $\Delta \xi A_x=\Delta \eta A_y=0.2$, $N_x=N_y=16$. We take the final time $t = 10$.
		
		\item \textbf{Randomized grid}.  The randomized grid is generated by a $41 \times 41$
		uniform grids $[-0.5,0.5]\times[-0.5,0.5]$ with 10\% magnitude grid spacing in a random direction. The final time is $t = 10$.
				 
		\item \textbf{Moving wavy grid}. The grid is given as 
		\begin{equation}\label{moving_grid}
		\begin{aligned}
		 x_{i,j}(\tau)&=x_{min}+\Delta \xi \left[ (i-1)+A_x\left(1+0.1 \times \sin{2\pi \omega \tau}\right) \sin{\frac{N_y(j-1)\Delta \eta}{L_y}} \right],\\
		 y_{i,j}(\tau)&=y_{min}+\Delta \eta \left[ (j-1)+A_y\left(1+0.1 \times \sin{2\pi \omega \tau}\right) \sin{\frac{N_x(i-1)\Delta \xi}{L_x}} \right],		
		\end{aligned}
		\end{equation}
		where,
		\begin{equation*}
		\begin{aligned}
		 &i=1,2,\cdots,I_{max},\qquad
		 j=1,2,\cdots,J_{max},\\
%		 &\Delta x_{0}=\frac{L_x}{I_{max}-1},\qquad 
%		 \Delta y_{0}=\frac{L_y}{J_{max}-1},\\	
		 &x_{min}=-\frac{L_{x}}{2},\qquad
		 y_{min}=-\frac{L_{y}}{2},
		\end{aligned}
		\end{equation*}
		with the specified parameters $I_{max}=J_{max}=21$, $L_x=L_y=1$, $\Delta \xi A_x=\Delta \eta A_y=0.05$, $N_x=N_y=4$, and the frequency of oscillation $\omega=1.0$. The final time is taken as $t = 10$.
		
		\item %{\color{red}
			\textbf{Stationary spherical grid}. The grid is expressed as
			\begin{equation}\label{cylinder_grid}
			\begin{aligned}
			x_{i,j}&=\left(r_1-\left(r_1-r_0\right)\Delta \xi(i-1)\right) \, \cos(\pi+\theta \, (1-2\Delta \eta(j-1))),\\
			y_{i,j}&=\left(r_2-\left(r_2-r_0\right)\Delta \xi(i-1)\right) \, \sin(\pi+\theta \, (1-2\Delta \eta(j-1))),		
			\end{aligned}
			\end{equation}
			where,
			\begin{equation*}
			\begin{aligned}
			&i=1,2,\cdots,I_{max},\qquad
			j=1,2,\cdots,J_{max},\\
			\end{aligned}
			\end{equation*}
			with $I_{max}=J_{max}=41$, $r_1=0.3$, $r_2=0.65$, $r_0=0.125$, $\theta=5\pi/12$ and $\Delta \xi =\Delta \eta =0.025$.
		 	And we test the problem with final time $t=0.1$.%}
		
	\end{itemize}

\begin{figure}[htbp]
		\centering	
		\subfigure[Stationary wavy grid]{
			\begin{minipage}[t]{0.4\linewidth}
				\centering
				\includegraphics[width=2in]{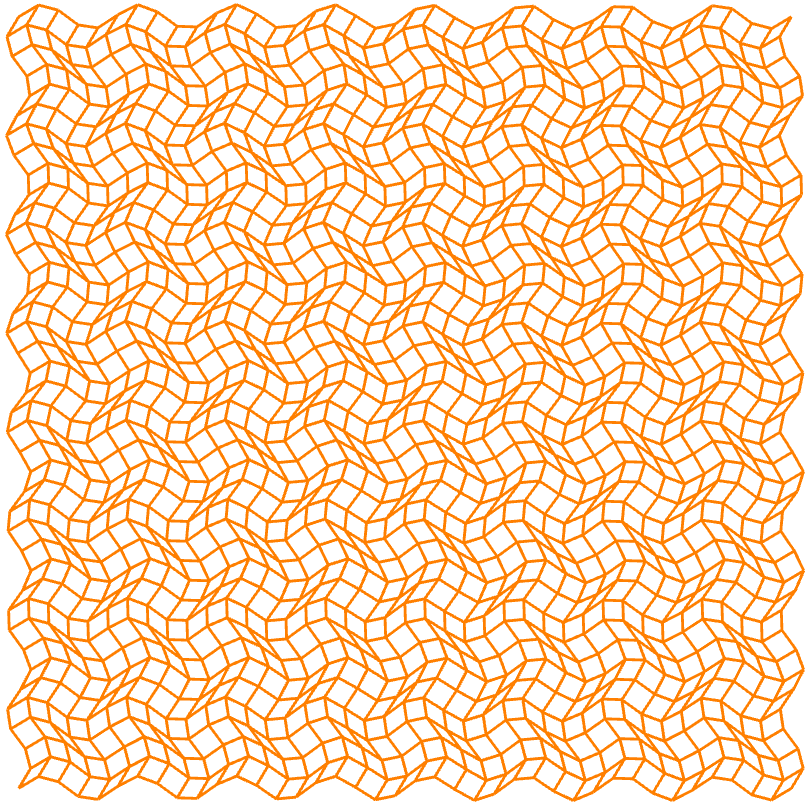}
				%\caption{fig1}
			\end{minipage}%
		}%	
		\subfigure[Randomized grid]{
			\begin{minipage}[t]{0.4\linewidth}
				\centering
				\includegraphics[width=2in]{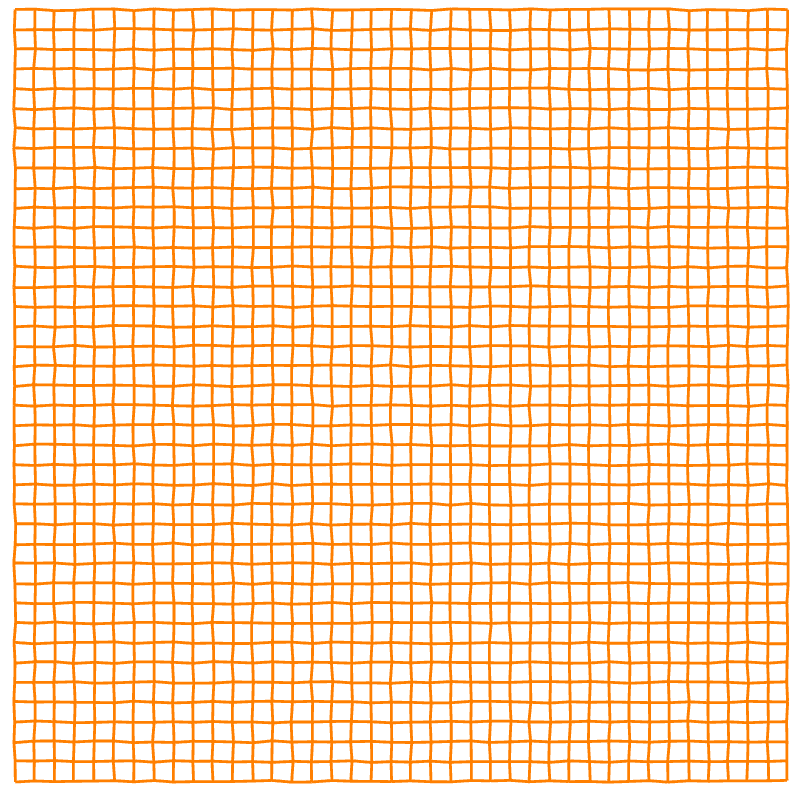}
				%\caption{fig1}
			\end{minipage}%
		} \\
		\subfigure[Moving wavy grid at $t=10$]{
			\begin{minipage}[t]{0.4\linewidth}
				\centering
				\includegraphics[width=2in]{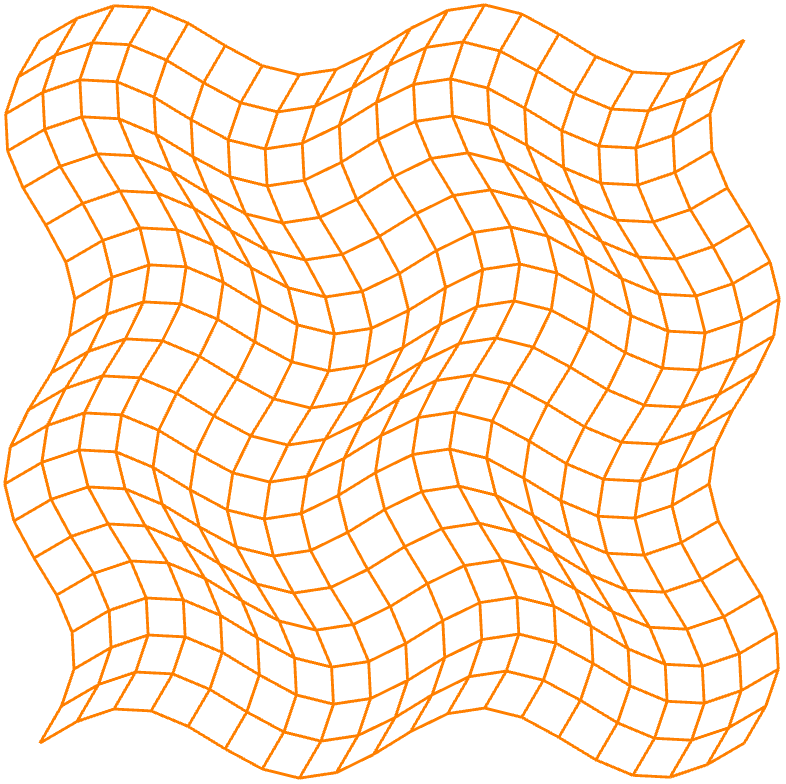}
				%\caption{fig2}
			\end{minipage}%
		}%		
	\subfigure[Stationary spherical grid]{
		\begin{minipage}[t]{0.4\linewidth}
			\centering
			\includegraphics[width=1.5in]{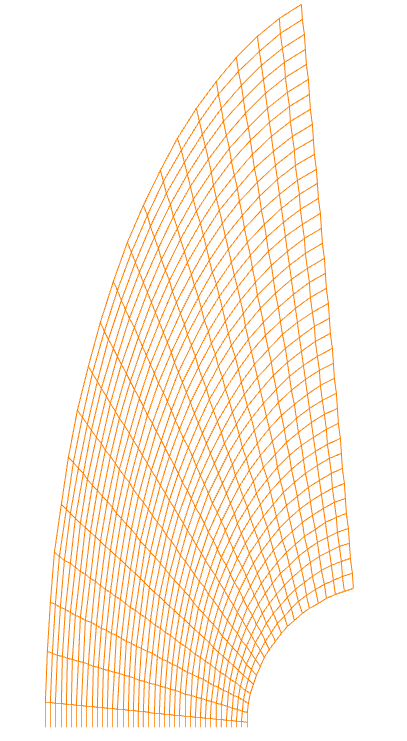}
			%\caption{fig2}
		\end{minipage}%
	}  
		\centering
		\caption{Example \ref{ex1}: The grids for the free-stream preserving test.}
		\label{grid1}
	\end{figure}

	The grids are shown in Fig \ref{grid1}. Time steps are taken as $\Delta \tau=0.05$, $\Delta \tau=0.05$, $\Delta \tau=0.01$ and $\Delta \tau=5.0 \times 10^{-4}$, respectively. The $L_\infty$ errors of $v$ at final time are shown in Table \ref{tab1}. We can see that PL-WENO have errors that are close to the machine zero. However, NPL-WENO have large errors on the level of $10^{-2}$. These demonstrate that PL-WENO can preserve the free-stream condition, while NPL-WENO can not.

%{\color{red}
%Moreover, we give an numerical example in stationary cylinder grid which is an analysis non-uniform mesh,
%and the grid is expressed as
%	\begin{equation}\label{cylinder_grid}
%	\begin{aligned}
%	x_{i,j}&=\left(r_1-\left(r_1-r_0\right)\Delta \xi(i-1)\right) \, \cos(\pi+\theta \, (1-2\Delta \eta(j-1))),\\
%	y_{i,j}&=\left(r_2-\left(r_2-r_0\right)\Delta \xi(i-1)\right) \, \sin(\pi+\theta \, (1-2\Delta \eta(j-1))),		
%	\end{aligned}
%	\end{equation}
%	where,
%	\begin{equation*}
%	\begin{aligned}
%	&i=1,2,\cdots,I_{max},\qquad
%	j=1,2,\cdots,J_{max},\\
%	\end{aligned}
%	\end{equation*}
%	with $I_{max}=J_{max}=41$, $r_1=0.3$, $r_2=0.65$, $r_0=0.125$, $\theta=5\pi/12$ and $\Delta \xi =\Delta \eta =0.025$.
%The grid are shown in Fig \ref{grid2}. Time step is taken as $\Delta \tau=5.0 \time 10^{-4}$. The $L_1$ and $L_\infty$ errors of $v$ at $t = 0.1$ are shown in Table \ref{tab2}.
% }

%	\begin{figure}[htbp]
%	\centering	
%	\subfigure[Stationary cylinder grid]{
%		\begin{minipage}[t]{0.3\linewidth}
%			\centering
%			\includegraphics[width=1.7in]{fig/2_4.png}
%			%\caption{fig2}
%		\end{minipage}%
%	}%	  
%	\centering
%	\caption{Example \ref{ex1}: The grids for the free-stream preserving test.}
%	\label{grid2}
%	\end{figure}

\begin{table}[htbp] 
	\caption{Example \ref{ex1}: $L_\infty$ errors of $v$ component in the free-stream preservation test.}
	\vspace{0.3cm}
	\centering
	\begin{small}
		\begin{tabular}{ccccc}
			\hline
			& Stationary wavy grid  & Randomized grid   &  Moving wavy grid &  Stationary spherical grid \\\hline
%			&  $L_\infty$ errors  &  $L_\infty$ errors  &  $L_\infty$ errors  &  $L_\infty$ errors \\ \hline
			PL-WENO    & 3.41E-14 & 8.21E-14  & 9.13E-14 & 6.32E-14  \\
			NPL-WENO & 6.85E-02 & 1.52E-02 & 2.16E-03 & 5.18E-07  \\			
			\hline
		\end{tabular}
	\end{small}
	\label{tab1}
\end{table}

%\begin{table}[htbp] 
%	\caption{Example \ref{ex1}: errors of $v$ component in the free-stream preservation test.}
%	\vspace{0.3cm}
%	\centering
%	\begin{small}
%	\begin{tabular}{ccccccccc}
%		\hline
%		 & \multicolumn{2}{c}{Stationary wavy grid}  & \multicolumn{2}{c}{Randomized grid}   &  \multicolumn{2}{c}{Moving wavy grid} &  \multicolumn{2}{c}{Stationary wavy grid} \\\cline{2-9}
%		 & $L_1$ errors &  $L_\infty$ errors  & $L_1$ errors &  $L_\infty$ errors  & $L_1$ errors &  $L_\infty$ errors  & $L_1$ errors &  $L_\infty$ errors \\ \hline
%		PL-WENO    & 8.12E-15 & 3.41E-14 & 7.57E-15 & 8.21E-14  & 2.32E-14 & 9.13E-14 & 5.75E-15 & 6.32E-14  \\
%		NPL-WENO & 1.78E-02 & 6.85E-02 & 2.95E-03 & 1.52E-02 & 8.03E-04 & 2.16E-03 & 6.44E-08 & 5.18E-07  \\			
%		\hline
%	\end{tabular}
%    \end{small}
%	\label{tab1}
%\end{table}

%\begin{table}[htbp] 
%	\caption{Example \ref{ex1}: errors of $v$ component in the free-stream preservation test.}
%	\vspace{0.3cm}
%	\centering
%	\begin{small}
%		\begin{tabular}{ccc}
%			\hline
%			%& \multicolumn{2}{c}{Stationary wavy grid}  \\\cline{2-3}
%			& $L_1$ errors &  $L_\infty$ errors \\ \hline
%			PL-WENO    & 5.75E-15 & 6.32E-14  \\
%			NPL-WENO & 6.44E-08 & 5.18E-07  \\			
%			\hline
%		\end{tabular}
%	\end{small}
%	\label{tab2}
%\end{table}

	\subsection{MHD: 2D field loop} \label{ex4}
	In this section we consider a 2D advection of a weakly magnetized field loop from \cite{gardiner2005unsplit}. The initial conditions are
	\begin{equation*}
	(\rho,u,v,w,p)(x,y,0)=(1,\sqrt 5\cos(\theta),\sqrt 5\sin(\theta),0,1),
	\end{equation*}
	with the advection angle of $\theta = \arctan (0.5)$. Magnetic field components are initialized by taking the curl of the magnetic potential $A$,
	\begin{equation*}
	A(x,y,0)=
	\left\{
	\begin{array}{ll}
	0.001(R-r), &\text{if} \ r\leq R,\\
	0, & \text{otherwise}, \\
	\end{array}
	\right.
	\end{equation*}
	with $r=\sqrt{x^2+y^2}$ and $R=0.3$. Periodic boundary conditions are used in both directions.
	
	We will test 2D field loop problem on wavy grid, randomized	grid and moving wavy grid. 
	In this problem, the wavy grid (\ref{wavy_grid}) and the moving wavy grid (\ref{moving_grid}) are formulated with the parameters 
	$$I_{max}=101, \quad J_{max}=51, \quad L_x=2, \quad L_y=1,  $$ 
	$$\Delta \xi A_x=0.03, \quad \Delta \eta A_y=0.06, \quad N_x=N_y=8,$$
	and the frequency of oscillation $\omega=1.0$.
	The randomized	grid is generated by a $101 \times 51$ uniform grids on $[-1,1]\times[-0.5,0.5]$ with $10\%$ magnitude grid spacing in a random direction. 

	The time step is chosen as $\Delta \tau= CFL \,\Delta \xi / \alpha$, with $CFL=0.1$, and $\alpha$ is the largest wave speed estimated in the curvilinear coordinates.
	
\begin{figure}[htbp]
		\centering	
		\subfigure[PL-WENO. Stationary wavy grid. ]{
			\begin{minipage}[t]{0.5\linewidth}
				\centering
				\includegraphics[width=3in]{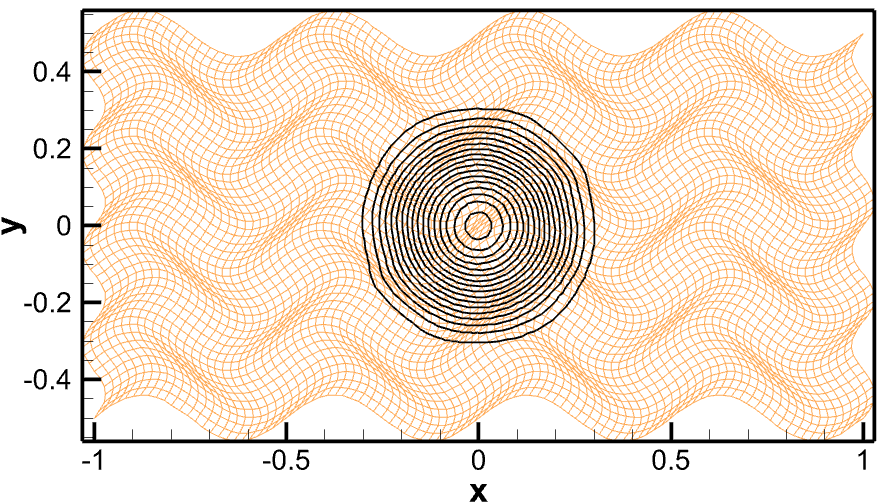}
				%\caption{fig1}
			\end{minipage}%
		}%
		\subfigure[NPL-WENO. Stationary wavy grid]{
			\begin{minipage}[t]{0.5\linewidth}
				\centering
				\includegraphics[width=3in]{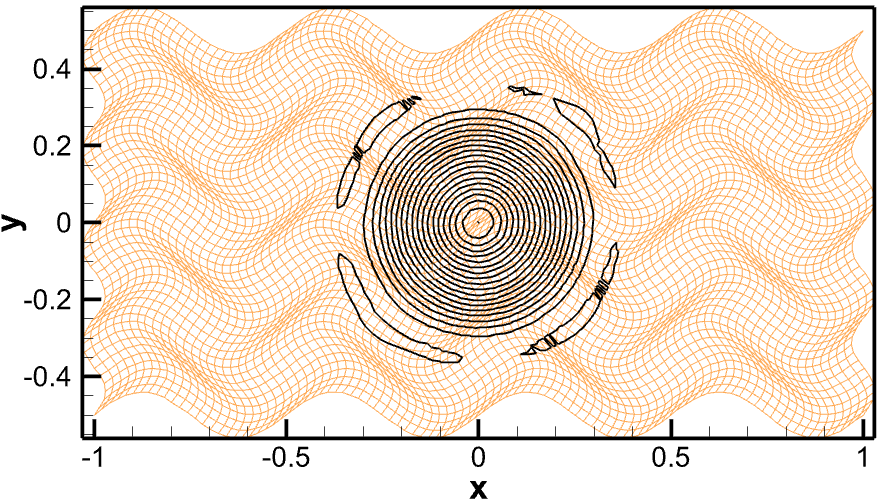}
				%\caption{fig2}
			\end{minipage}%
		}\\
	\subfigure[PL-WENO. Randomized grid.]{
		\begin{minipage}[t]{0.5\linewidth}
			\centering
			\includegraphics[width=3in]{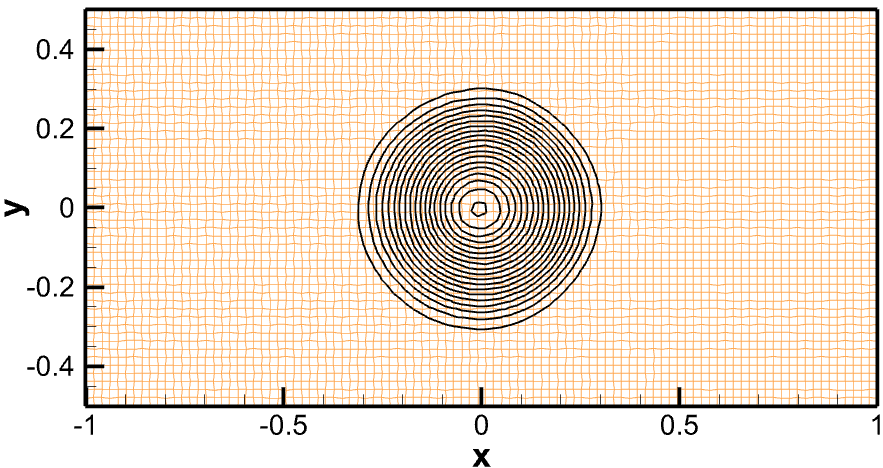}
			%\caption{fig1}
		\end{minipage}%
	}%
	\subfigure[NPL-WENO. Randomized grid.]{
		\begin{minipage}[t]{0.5\linewidth}
			\centering
			\includegraphics[width=3in]{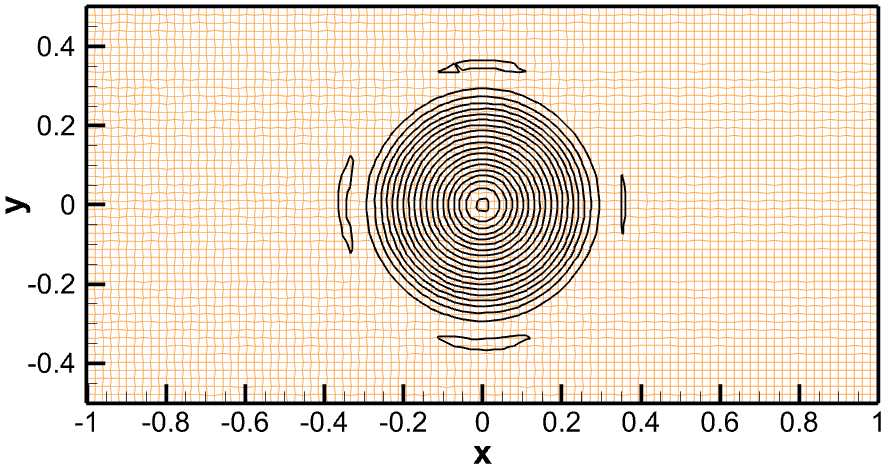}
			%\caption{fig2}
		\end{minipage}%
	}\\%	%	
	\subfigure[PL-WENO. Moving wavy grid.]{
	\begin{minipage}[t]{0.5\linewidth}
		\centering
		\includegraphics[width=3in]{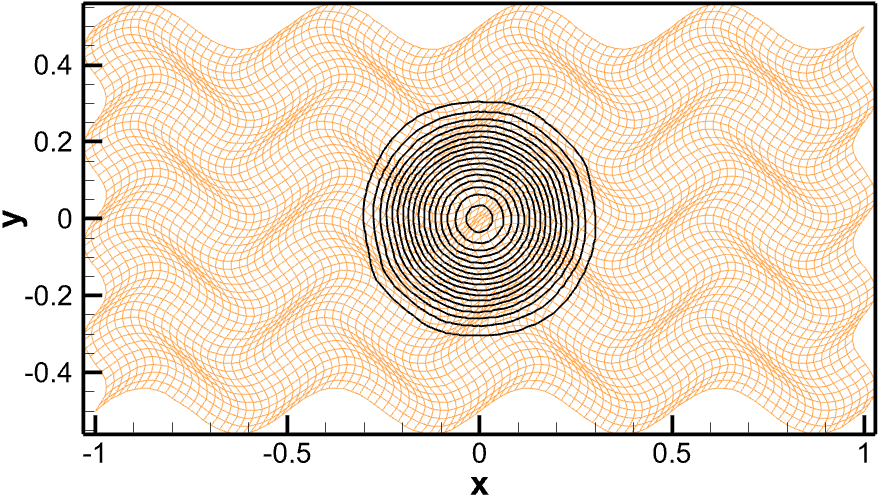}
		%\caption{fig1}
	\end{minipage}%
	}%
\subfigure[NPL-WENO. Moving wavy grid.]{
	\begin{minipage}[t]{0.5\linewidth}
		\centering
		\includegraphics[width=3in]{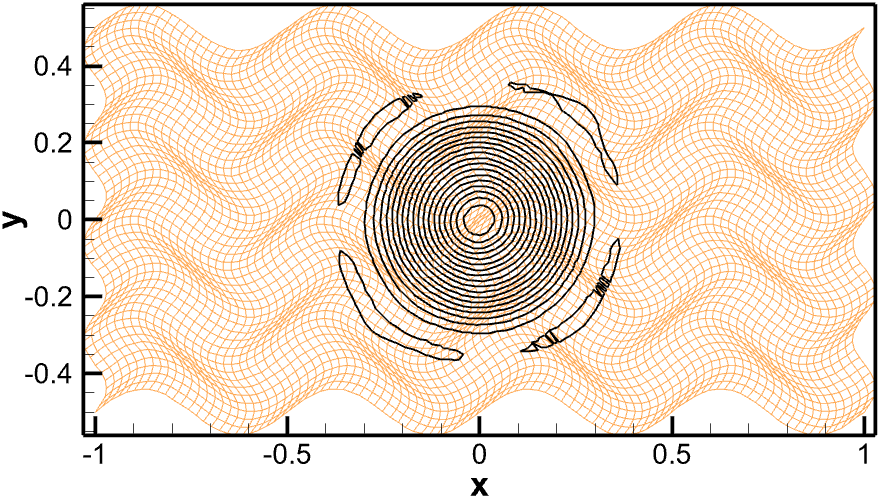}
		%\caption{fig2}
	\end{minipage}%
}%	
\centering
		\caption{Example \ref{ex4}: $A$ of the 2D field loop. 20 contour lines between $-2.16\times 10^{-6}$ and $2.7\times 10^{-4}$ are shown.}
		\label{example4_1}
	\end{figure}

\begin{figure}[htbp]
		\centering
		\subfigure[PL-WENO. Stationary wavy grid.]{
			\begin{minipage}[t]{0.5\linewidth}
				\centering
				\includegraphics[width=3in]{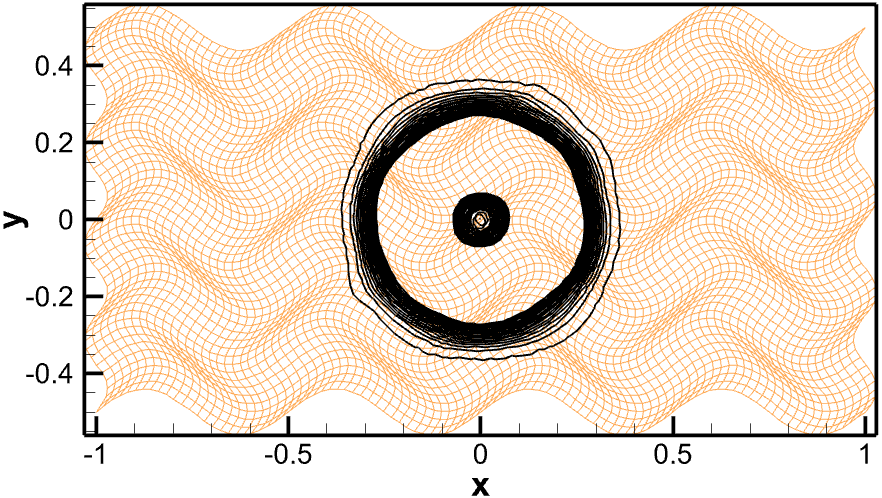}
				%\caption{fig1}
			\end{minipage}%
		}%
		\subfigure[NPL-WENO. Stationary wavy grid.]{
			\begin{minipage}[t]{0.5\linewidth}
				\centering
				\includegraphics[width=3in]{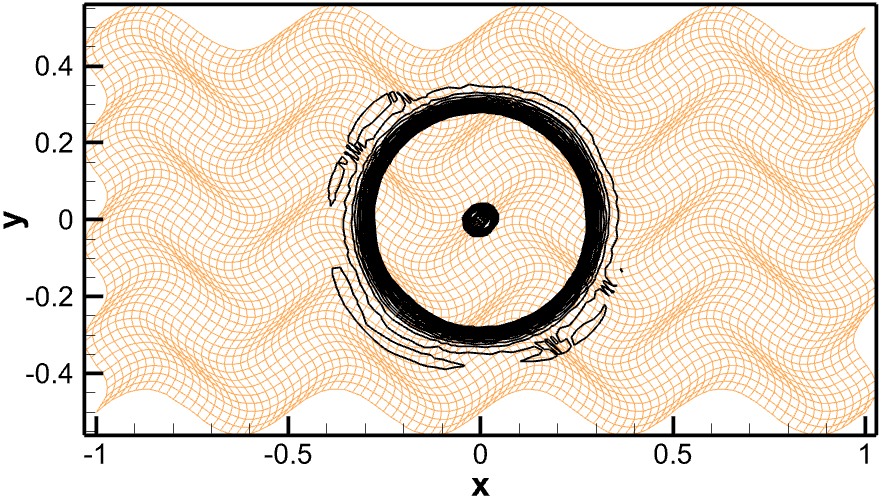}
				%\caption{fig2}
			\end{minipage}%
		} \\ %
		\subfigure[PL-WENO. Randomized grid.]{
			\begin{minipage}[t]{0.5\linewidth}
				\centering
				\includegraphics[width=3in]{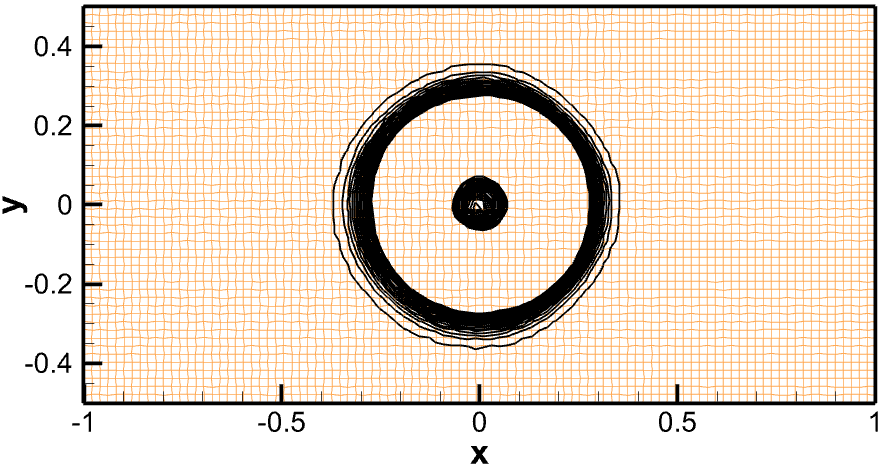}
				%\caption{fig1}
			\end{minipage}%
		}%
		\subfigure[NPL-WENO. Randomized grid.]{
			\begin{minipage}[t]{0.5\linewidth}
				\centering
				\includegraphics[width=3in]{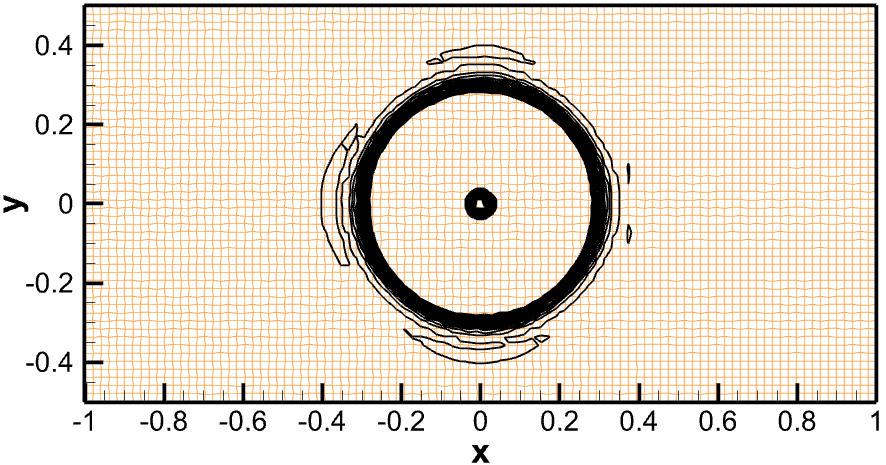}
				%\caption{fig2}
			\end{minipage}%
		} \\%
		\subfigure[PL-WENO. Moving wavy grid.]{
			\begin{minipage}[t]{0.5\linewidth}
				\centering
				\includegraphics[width=3in]{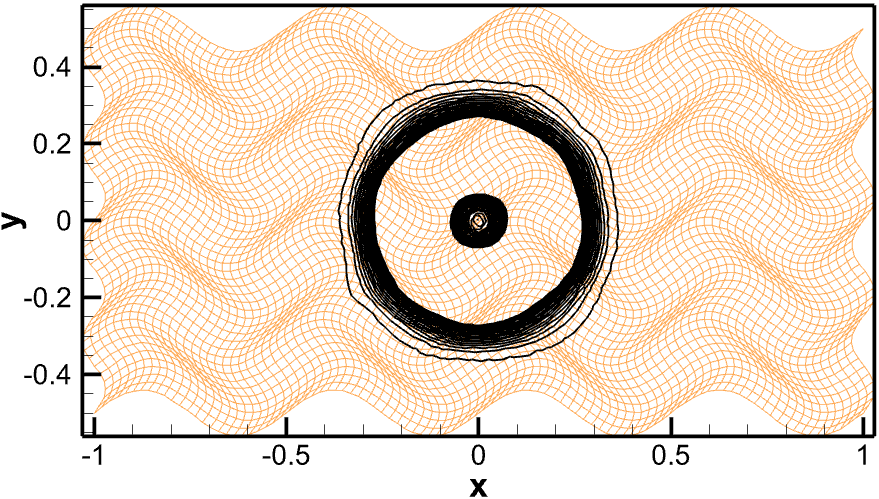}
				%\caption{fig1}
			\end{minipage}%
		}%
		\subfigure[NPL-WENO. Moving wavy grid.]{
			\begin{minipage}[t]{0.5\linewidth}
				\centering
				\includegraphics[width=3in]{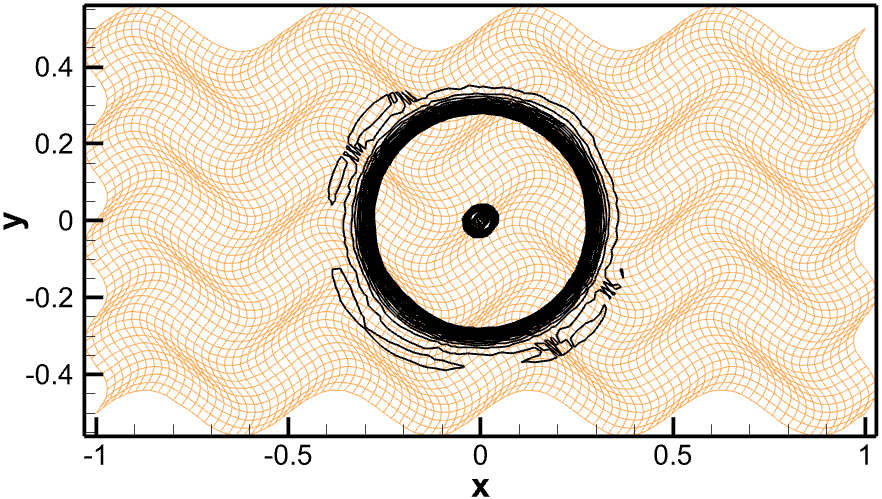}
				%\caption{fig2}
			\end{minipage}%
		}%
		\centering
		\caption{Example \ref{ex4}: $\|\textbf{B}\|^2$ of the 2D field loop. 20 contour lines between $3\times 10^{-9}$ and $5.2\times 10^{-7}$ are shown}
		\label{example4_2}
	\end{figure}

	We present contour plots of the potential $A$ and $\|\textbf{B}\|^2=B_1 ^2+B_2 ^2+B_3^2$ at $t=2$ in Fig. \ref{example4_1} and Fig. \ref{example4_2}. It is observed that both schemes can maintains the circular symmetry of the loop as expected, but there are small irrational disturbances founded in NPL-WENO solution.
	
	\subsection{MHD: 2D rotor problem}\label{ex5}
\begin{figure}[htbp]
		\centering	
		\subfigure[PL-WENO]{
			\begin{minipage}[t]{0.5\linewidth}
				\centering
				\includegraphics[width=2.5in]{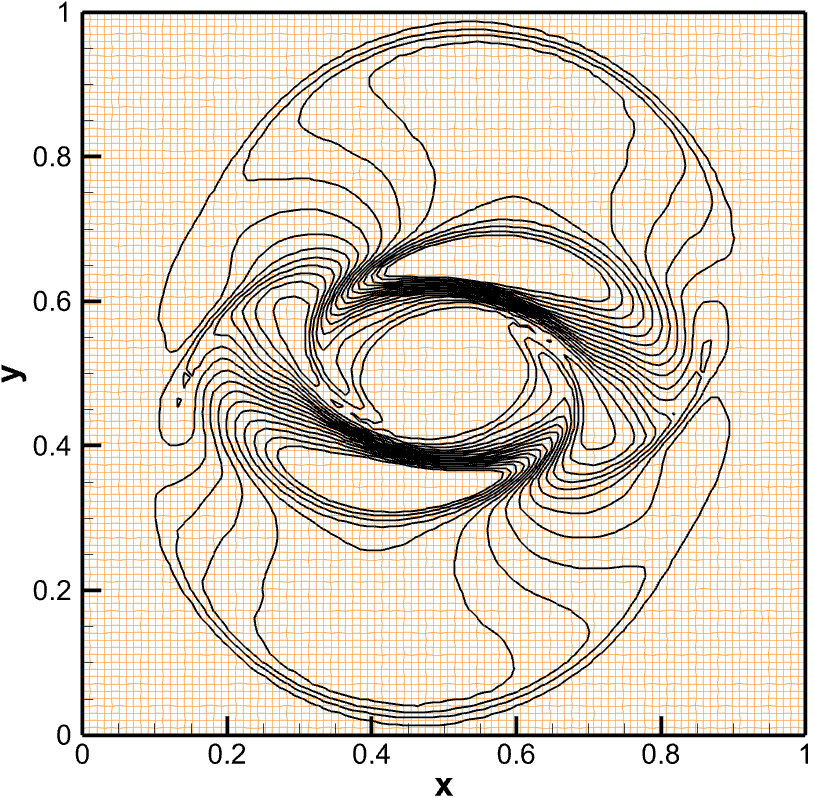}
				%\caption{fig1}
			\end{minipage}%
		}%
		\subfigure[NPL-WENO]{
			\begin{minipage}[t]{0.5\linewidth}
				\centering
				\includegraphics[width=2.5in]{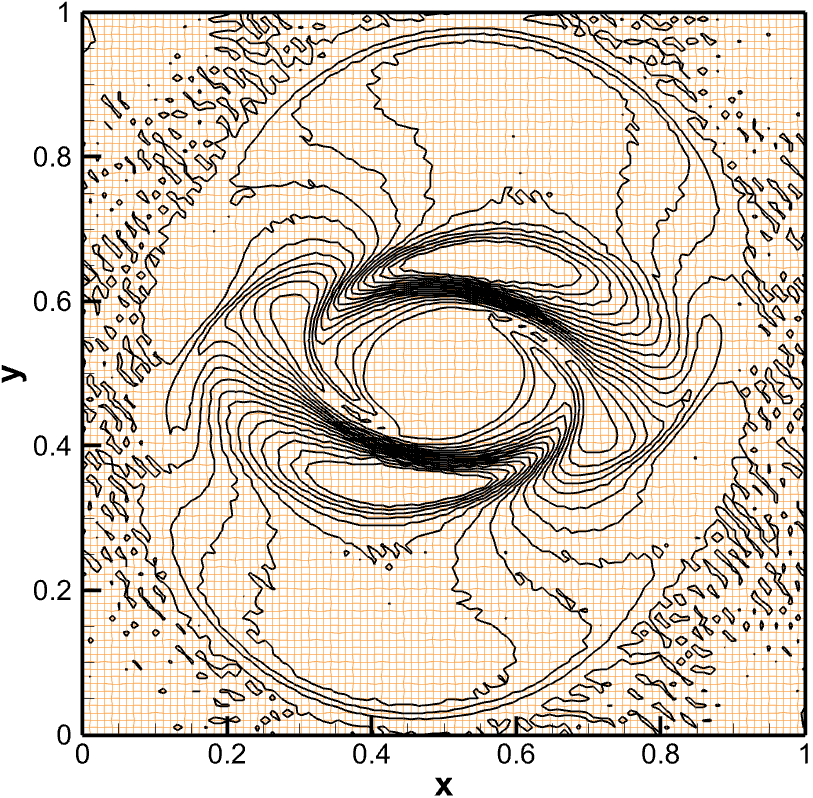}
				%\caption{fig2}
			\end{minipage}%
		}%	
		\centering
		\caption{Example \ref{ex5}: $\|\textbf{B}\|^2$ of the 2D rotor problem. 16 contour lines between $0.1$ and $1.0$ are shown.}
		\label{example4.5}
	\end{figure}

	The initial condition is given as
	\begin{equation*}
	\left(\rho,u,v\right)
	=
	\left\{
	\begin{array}{ll}
	\left(
	10,-\left(y-0.5\right)/r_0,\left(x-0.5\right)/r_0
	\right),\quad 
	&\text{if} \ r \leq r_0, \\
	\left(
	1+9f(r),-f(r)(y-0.5)/r,f(r)(x-0.5)/r
	\right),\quad 
	&\text{if} \  r_0 \le r \leq r_1,\\
	\left(
	1,0,0
	\right),\quad 
	&\text{if} \   r \ge r_1,\\
	\end{array}
	\right.
	\end{equation*}
	and
	\begin{equation*}
	w=0,\quad B_1=2.5 /\sqrt{4\pi}, \quad B_2=0, \quad B_3=0, 
	\quad p=0.5, \quad A=2.5/\sqrt{4\pi} y,
	\end{equation*}
	where $r=\sqrt{(x-0.5)^2+(y-0.5)^2}$, $r_0=0.1$, $r_1=0.115$ and
	$f(r)=(r_1-r)/(r_1-r_0)$.
	Here we use the same initial condition of the second rotor problem test in \cite{toth2000b}. The problem is solved on a randomized grid generated by a $101 \times 101$ uniform grids with $5\%$ magnitude grid spacing in a random direction.	
	Same method of time step chosen as in 2D field loop is used here. 
	The contour plots of $\|\textbf{B}\|^2=B_1 ^2+B_2 ^2+ B_3^2$ at $t=0.295$ are presented in Fig. \ref{example4.5}. Results of NPL-WENO show large numerical errors owing to the grid distortions, but PL-WENO can ignore the discontinuity in the derivatives of the randomized grid to some extent from the result.
	
\subsection{MHD: 2D blast wave} \label{ex3}
Next we consider the blast wave problem. In this test strong shocks interact with a low-$\beta$ background, which could potentially cause negative density or pressure in numerical simulations. The problem has been commonly used to test the positivity-preserving capabilities of numerical methods for MHD equation \cite{balsara2009efficient, balsara1999staggered}. The initial conditions contain a constant density, velocity and magnetic field
\begin{equation*}
(\rho,u,v,w,B_1,B_2,B_3)(x,y,0)=(1,0,0,0,50\sqrt{2\pi},50\sqrt{2\pi},0),
\end{equation*}
and a piecewise defined pressure
\begin{equation*}
p(x,y,0)=
\left\{
\begin{array}{ll}
1000, & \text{if} \ r\leq 0.1,\\
0.1, & \text{otherwise}. \\
\end{array}
\right.
\end{equation*}
Similarly, we will test 2D blast wave problem on wavy grid, randomized grid and moving wavy grid. 
%In this problem, the wavy grid is mapped from the computational domain $(\xi,\eta)\in [0,1]\times[0,1]$ formulated as (\ref{wavy_grid}), and 
For both the wavy grid (\ref{wavy_grid}) and the moving wavy grid (\ref{moving_grid}), the parameters are taken as 
$$I_{max}=J_{max}=101, \quad L_x=L_y=1, \quad  N_x=N_y=4, \quad \Delta \xi A_x=0.03=\Delta \eta A_y=0.03, $$
and the frequency of oscillation $\omega=1.0$. The randomized grid is generated by a $101 \times 101$ uniform grids on $[0,1]\times[0,1]$ with $1\%$ magnitude grid spacing in a random direction. 
Same method of time step chosen as in 2D field loop is used here.

\begin{figure}[htbp]
	\centering
	\subfigure[PL-WENO. Stationary wavy grid. $T=0.0025$.]{
		\begin{minipage}[t]{0.5\linewidth}
			\centering
			\includegraphics[width=2.5in]{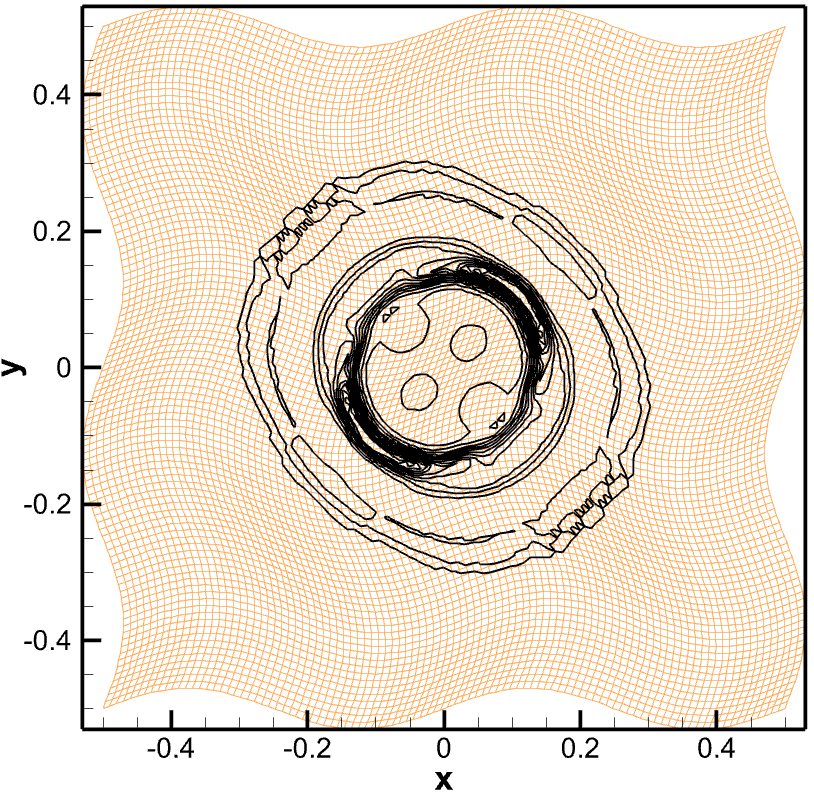}
			%\caption{fig1}
		\end{minipage}%
	}%
	\subfigure[NPL-WENO. Stationary wavy grid. $T=0.0025$.]{
		\begin{minipage}[t]{0.5\linewidth}
			\centering
			\includegraphics[width=2.5in]{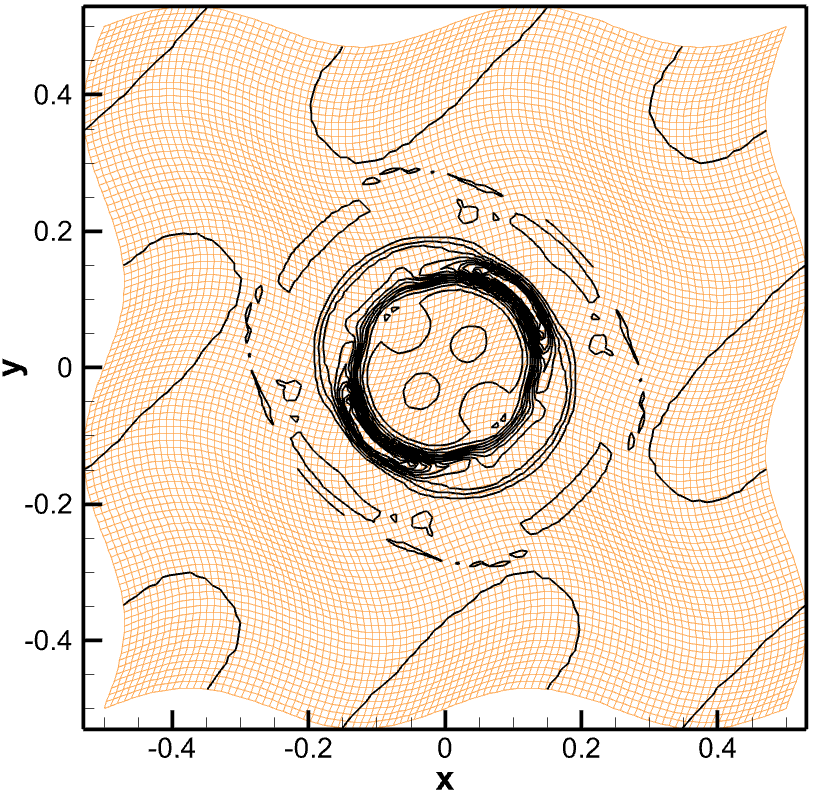}
			%\caption{fig2}
		\end{minipage}%
	}\\ %
	\subfigure[PL-WENO. Randomized grid. $T=0.0015$.]{
		\begin{minipage}[t]{0.5\linewidth}
			\centering
			\includegraphics[width=2.5in]{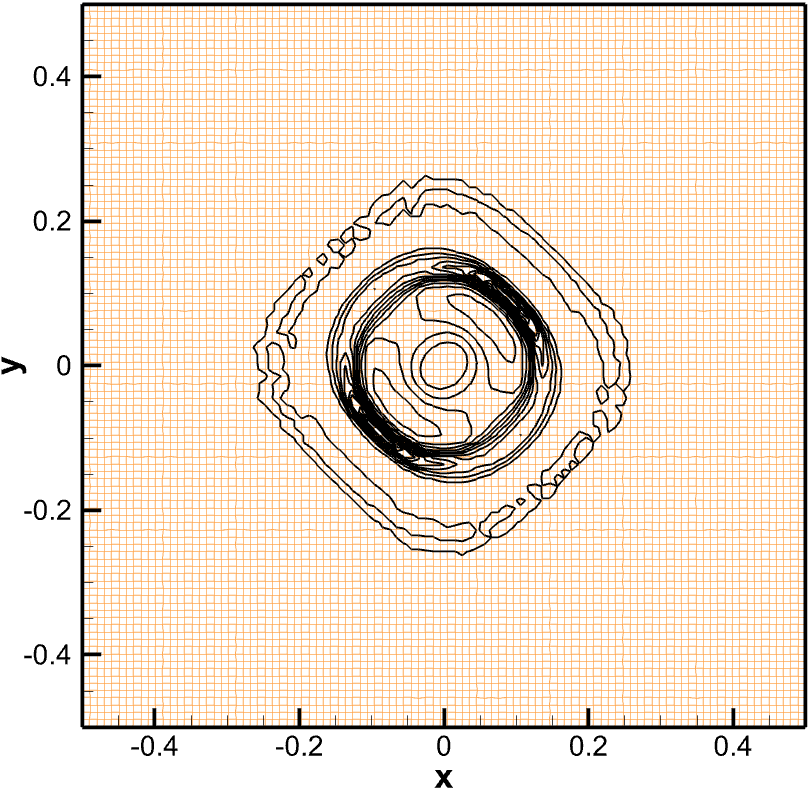}
			%\caption{fig1}
		\end{minipage}%
	}%
	\subfigure[NPL-WENO. Randomized grid. $T=0.0015$.]{
		\begin{minipage}[t]{0.5\linewidth}
			\centering
			\includegraphics[width=2.5in]{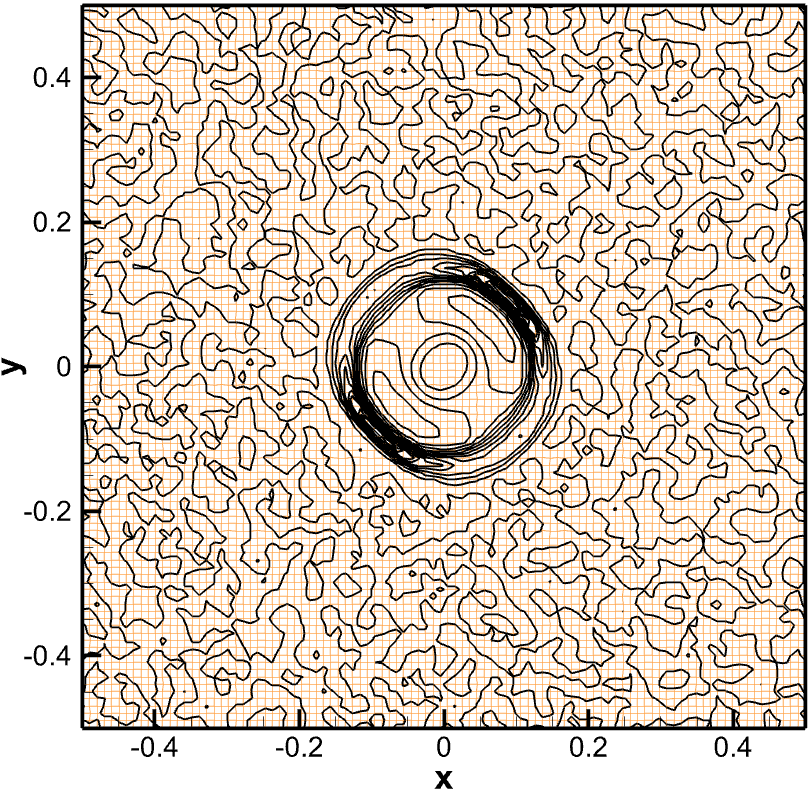}
			%\caption{fig2}
		\end{minipage}%
	} \\ %
	\subfigure[PL-WENO. Moving wavy grid. $T=0.0013$.]{
		\begin{minipage}[t]{0.5\linewidth}
			\centering
			\includegraphics[width=2.5in]{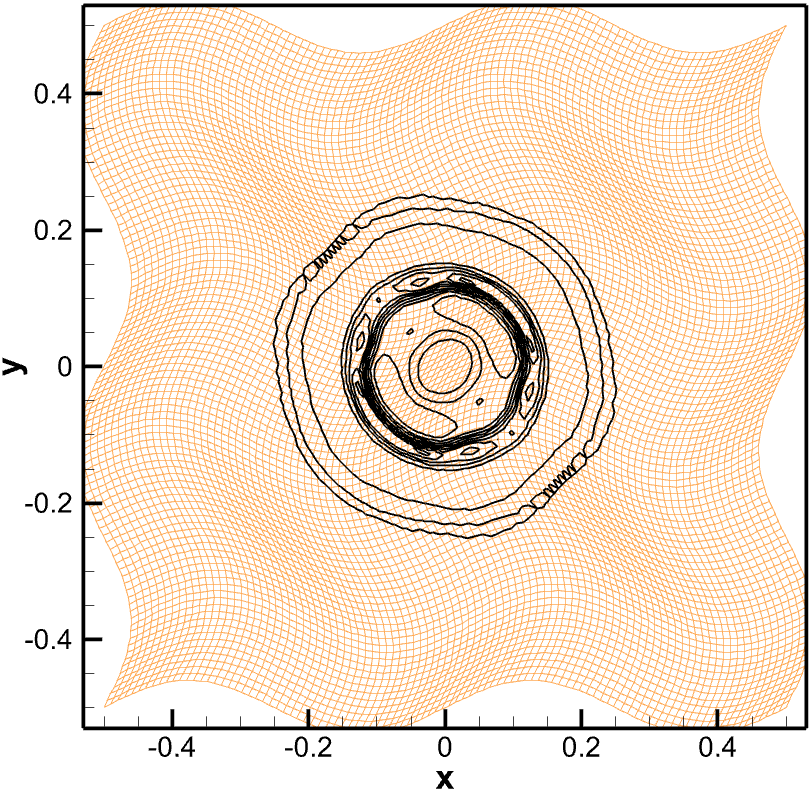}
			%\caption{fig1}
		\end{minipage}%
	}%
	\subfigure[NPL-WENO. Moving wavy grid. $T=0.0013$. ]{
		\begin{minipage}[t]{0.5\linewidth}
			\centering
			\includegraphics[width=2.5in]{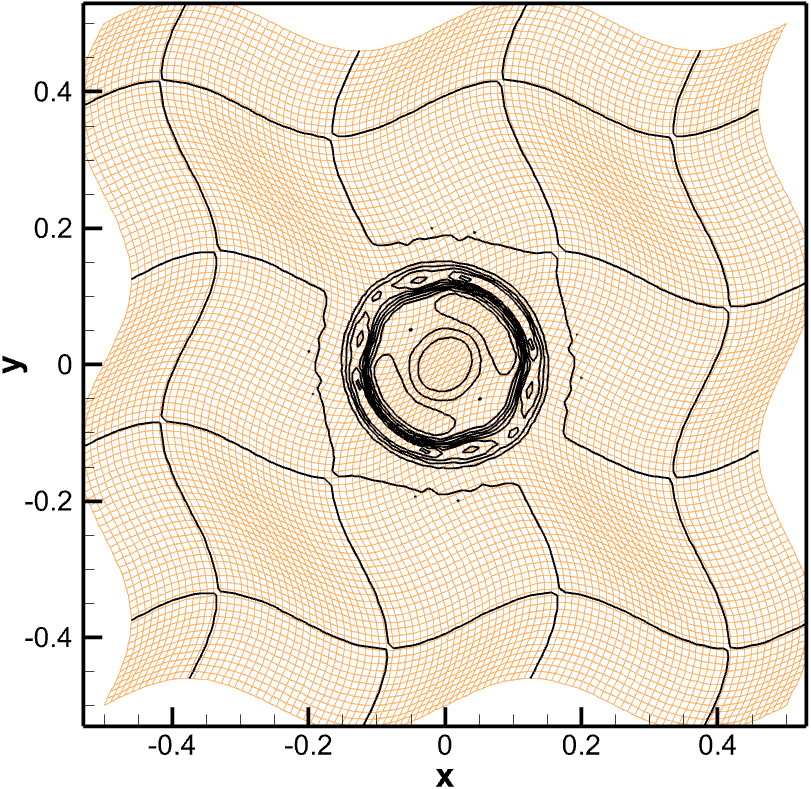}
			%\caption{fig2}
		\end{minipage}%
	}%
	\centering
	\caption{Example \ref{ex3}: $\rho$ of the 2D blast wave. In (a)-(b), we use 15 contour lines between $0.5$ and $1.9$ are shown; in (c)-(f), we use 19 contour lines between $0.5$ and $1.9$ are shown. }
	\label{example9}
\end{figure}

%{\color{red} 
	PL-WENO can successfully run to $T=0.01$, keeping density and pressure positive all the time in the computational domain for all three types of grids.%}
However, NPL-WENO would produce negative pressures at around  $T=0.002898$, $T=0.002257$, $T=0.001486$ on the three types of grids, respectively. 
This indicates that missing free-stream-preserving property can lead large error.
Moreover, in Fig.\ref{example9}, we present contour plots of $\rho$ at $T=0.0025$ on wavy grid, $T=0.0015$ on randomized grid and $T=0.0013$ on moving wavy grid. It's obvious that NPL-WENO  would introduce extra errors and can not maintain the structure of solution on the three types of grids, while PL-WENO can.

%{\color{red}
\subsection{MHD: Bow shock flow} \label{ex41}
In this test,
we simulate the supersonic flow past a cylinder,
to examine the performance of our scheme when applied to problems involving curved physical boundaries.
The computational grid is generated by (\ref{cylinder_grid})
%The computational domain is determined by the mapping
%\begin{equation*}\label{cylinder_grid2}
%	\begin{aligned}
%		x&=\left(r_1-\left(r_1-r_0\right)\xi\right)\cos(\pi+\left(1-2\eta\right)\theta),\\
%		y&=\left(r_2-\left(r_2-r_0\right)\xi\right)\sin(\pi+\left(1-2\eta\right)\theta),	
%	\end{aligned}
%\end{equation*}
with $r_1=0.3$, $r_2=0.65$, $r_0=0.125$, $\theta=5\pi/12$ and  $\Delta \xi =\Delta \eta =\frac{1}{150}$. 
Assuming $r=\sqrt{x^{2}+y^{2}}$, $\delta r=0.125$, 
the initial condition is given as 
\begin{equation*}
\begin{aligned}
(\rho,u,v,w,p)=(1,2,0,0,0.2),\\
\mathbf{B}=\left\{
\begin{array}{ll}{
	\left(B_{1}, B_{2}, 0\right)^{T}}, & {\text { if } r \leq r_{0}+\delta r}, \\ 
	{(0.1,0,0)^{T}}, & {\text { otherwise, }}
	\end{array}
	\right.
\end{aligned}
\end{equation*}
where,

\begin{align*}
B_{1}&=0.1 \frac{\pi y^{2}}{2 \delta r \cdot r} \cos \left(\frac{\pi\left(r-r_{0}\right)}{2 \delta r}\right)+0.1 \sin \left(\frac{\pi\left(r-r_{0}\right)}{2 \delta r}\right),\\
B_{2}&=-0.1 \frac{\pi x y}{2 \delta r \cdot r} \cos \left(\frac{\pi\left(r-r_{0}\right)}{2 \delta r}\right).
\end{align*}
The corresponding initial magnetic potential is
\begin{equation}
\begin{aligned}
A=\left\{
\begin{array}{ll}
0.1 y \sin \left(\frac{\pi\left(r-r_{0}\right)}{2 \delta r}\right), &{\text { if } r \leq r_{0}+\delta r}, \\
0.1 y,  &\text{otherwise}.
\end{array}
\right.
\end{aligned}
\end{equation}
Here we use a uniform grid of size $150 \times 150$ in the domain $(\xi,\eta)$, and the boundary condition is identical to the bow shock flow test in \cite{christlieb2018high}.

\begin{figure}[htbp]
	\centering
	\subfigure[PL-WENO.]{
		\begin{minipage}[t]{0.5\linewidth}
			\centering
			\includegraphics[width=2.5in]{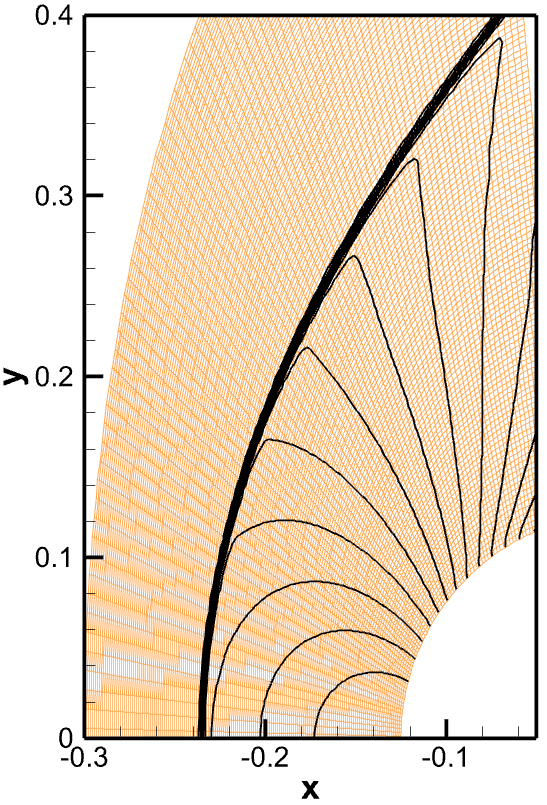}
			%\caption{fig1}
		\end{minipage}%
	}%
	\subfigure[NPL-WENO.]{
		\begin{minipage}[t]{0.5\linewidth}
			\centering
			\includegraphics[width=2.5in]{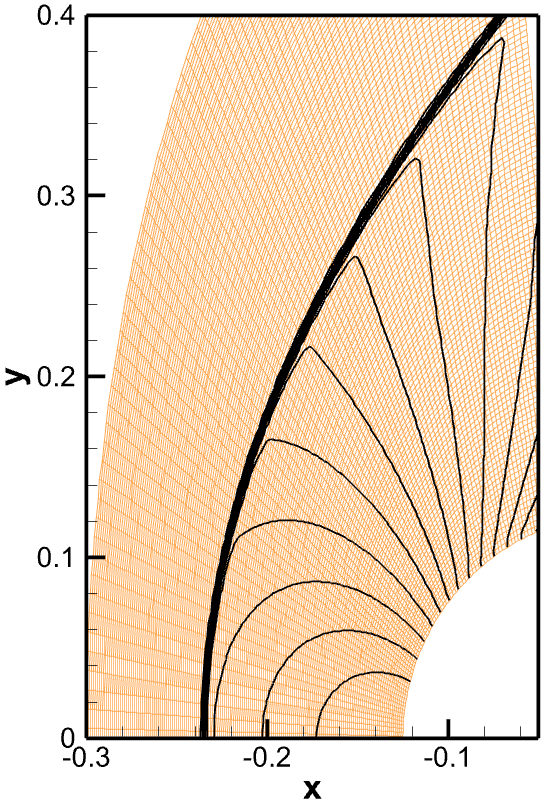}
			%\caption{fig2}
		\end{minipage}%
	}
	\centering
	\caption{Example \ref{ex41}: pressure of the bow shock flow at $t=6$. In (a)-(b), we use 16 contour lines between $0.4$ and $3.4$ are shown. }
	\label{example12}
\end{figure}

\begin{figure}[htbp]
	\centering
	\subfigure[PL-WENO.]{
		\begin{minipage}[t]{0.5\linewidth}
			\centering
			\includegraphics[width=2.5in]{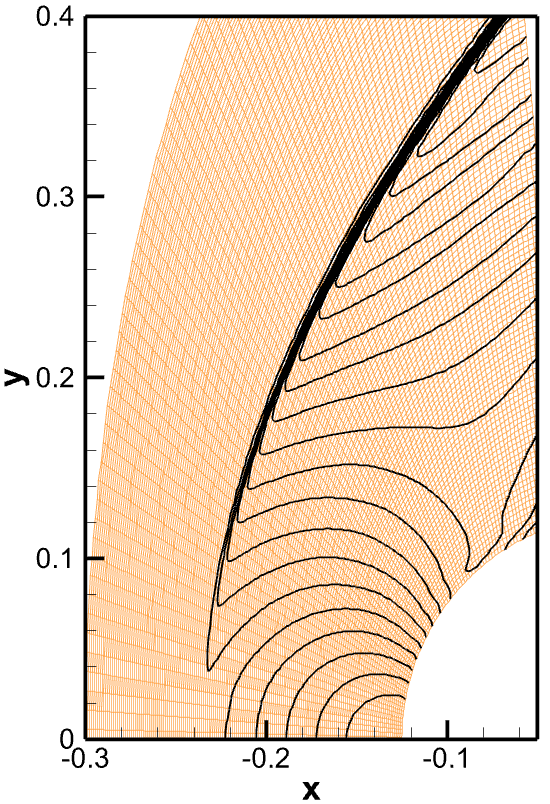}
			%\caption{fig1}
		\end{minipage}%
	}%
	\subfigure[NPL-WENO.]{
		\begin{minipage}[t]{0.5\linewidth}
			\centering
			\includegraphics[width=2.5in]{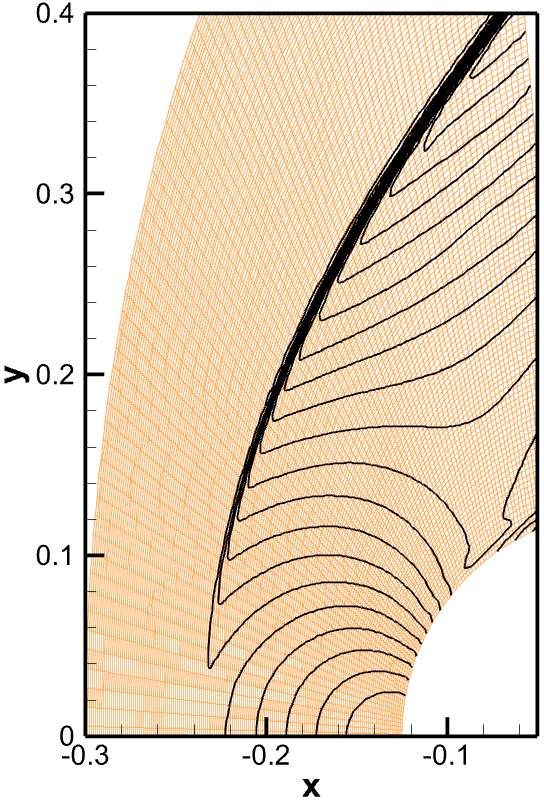}
			%\caption{fig2}
		\end{minipage}%
	}
	\centering
	\caption{Example \ref{ex41}: $\|\textbf{B}\|^2$ of the bow shock flow at $t=6$. In (a)-(b), we use 20 contour lines between $0.002$ and $0.034$ are shown. }
	\label{example13}
\end{figure}

The contour plots of pressure $p$ and norm of the magnetic field $\|\textbf{B}\|^2$ at $t=6$ are presented in Fig. \ref{example12} and Fig. \ref{example13}, respectively. 
We can see that both two schemes can simulate the bow shock flow well.
%but small irrational disturbances from norm of the magnetic field are founded in NPL-WENO solution.

%}

\section{Concluding remarks}	

In this paper, we have discussed the performance of finite difference WENO scheme on solving the ideal MHD equation on generalized meshes. 
A constrained transport framework is used to control the divergence error of the magnetic field. In this framework, a magnetic vector potential is evolved by a H-J equation, and set the magnetic field to be the (discrete) curl of the magnetic potential after each time step.
In order to preserve the free-stream condition, we employ an alternative flux formulation of the finite difference WENO scheme to the ideal MHD equations \cite{christlieb2018high}. 
%{\color{red}
	Moreover, a novel method is proposed to solve the H-J equations, which combines the standard WENO approximation on lattice meshes with a general Lax-Friedrichs type monotone Hamiltonian. 
In particular, the metric terms in the method are carefully defined such that they can be canceled exactly, and further ensure the scheme could maintain the linear function solutions of H-J equations on both stationary and dynamically generalized coordinate systems.%}
The free-stream preservation property of ideal MHD equation for the new schemes have been investigated both theoretically and numerically. 
Several numerical benchmark problems are used to confirm that the standard finite difference WENO schemes have a rather large error arising from the metric terms on randomized and moving grids, while the proposed WENO schemes can simulate the solution well for all grids.

\appendix
\section{WENO method}
For the sake of completeness, we give the formula of WENO interpolation for a hyperbolic system and WENO approximation for a H-J equation. Here, we only show the 1D case with fifth order accuracy. For 2D case, it can be obtained in a dimension-by-dimension fashion, in which we fix one index and do the process in the other direction.

\subsection{WENO interpolation for a system}
\label{sec:WENO}

Given the point values $\{\textbf{q}_i\}$, we will approximate the values at half points $\textbf{q}_{i+1/2}^{\pm}$ :

\begin{enumerate}
	\item Compute an average state $\textbf{q}_{i+1/2}$, using either the simple arithmetic mean or a Roe average.  Construct the right and left eigenvectors of the Jacobian $\partial \widetilde{\textbf{f}} / \partial \textbf{q}$,
	and denote their matrices by
	$$R_{i+1/2}=R(\textbf{q}_{i+1/2}) ,\quad R^{-1}_{i+1/2}=R^{-1}(\textbf{q}_{i+1/2}).$$
	
	\item Project the conserved quantities $\textbf{q}$ to the local characteristic variables $\textbf{v}$,
	\begin{equation*}
	\textbf{v}_{j} = R^{-1}_{i+1/2} \textbf{q}_{j}, \quad \text{for } j = i-2, \ldots, i+3.
	\end{equation*}
	
	\item Perform a scalar WENO interpolation on each component of the characteristic variable $\textbf{v}_{j}$ to obtain the corresponding component of $\textbf{v}_{i+1/2}^{\pm}$. 
	Here, the procedure of a fifth-order WENO interpolation to obtain the $k$-th component $v^{-}_{k,i+1/2}$ is described:
	\begin{enumerate}
		\item Choose one big stencil as $S = \{x_{i-2},\ldots,x_{i+2}\}$, and three small stencils as $S^{(r)} _{i}= \{x_{i-r},\ldots, x_{i-r+2}\}$, $r=0,1,2$.  
		On each small stencil,  the standard interpolation gives
			\begin{align*}
			& v^{(0)}_{k,i+1/2} = \frac{3}{8}v_{k,i} +\frac{3}{4}v_{k,i+1} -\frac{1}{8}v_{k,i+2},\\
			& v^{(1)}_{k,i+1/2} = -\frac{1}{8}v_{k,i-1} +\frac{3}{4}v_{k,i} +\frac{3}{8}v_{k,i+1},\\
			& v^{(2)}_{k,i+1/2} = \frac{3}{8}v_{k,i-2} -\frac{5}{4}v_{k,i-1} +\frac{15}{8}v_{k,i}, \\
			& v^\text{big}_{k,i+1/2} = d_{0} v^{(0)}_{k,i+1/2} +d_{1} v^{(1)}_{k,i+1/2} +d_{2} v^{(2)}_{k,i+1/2},
			\end{align*}
		with the linear weights being $d_{0}={5}/{16}$, $d_{1}={5}/{8}$ and $d_{2}={1}/{16}$.
		%$$d_{0}=\frac{5}{16}, \quad d_{1}=\frac{5}{8} , \quad d_{2}=\frac{1}{16}.$$
		
		\item Compute nonlinear weights $\omega_r$ from the linear weights $d_r$, 
		\begin{equation*}
		\omega_{r}=\frac{\alpha_{r}}{\alpha_{0}+\alpha_{1}+\alpha_{2}}, \quad \alpha_{r}=\frac{d_{r}}{(IS_{r}+\epsilon)^2}, \quad r=0,1,2,
		\end{equation*}
		where $\epsilon = 10^{-6}$ is used to avoid division by zero, and the smoothness indicators are given by
			\begin{align*}
			&IS_{0} = \frac{13}{12}\left( v_{k,i}-2v_{k,i+1}+v_{k,i+2}\right)^2 +\frac{1}{4}\left( 3v_{k,i}-4v_{k,i+1}+v_{k,i+2}\right)^2,\\
			&IS_{1} = \frac{13}{12}\left( v_{k,i-1}-2v_{k,i}+v_{k,i+1}\right)^2 +\frac{1}{4}\left( v_{k,i}-v_{k,i+1}\right)^2,\\
			&IS_{2} = \frac{13}{12}\left( v_{k,i-2}-2v_{k,i-1}+v_{k,i}\right)^2 +\frac{1}{4}\left( v_{k,i-2}-4v_{k,i-1}+3v_{k,i}\right)^2.
			\end{align*}

		\item The WENO interpolation is defined by
		$$v^{-}_{k,i+1/2}=\sum_{r=0}^{2}\omega_{r}v^{(r)}_{k,i+1/2}.$$
	\end{enumerate}
	\item Finally, project $\textbf{v}_{i+1/2}^{\pm}$ back to the conserved quantities,
	\begin{equation*}
	\textbf{q}_{i+1/2}^{\pm} = R_{i+1/2} \textbf{v}_{i+1/2}^{\pm}.
	\end{equation*}
\end{enumerate}
	
	Note that the process to obtain $\textbf{v}^{+}_{i+1/2}$ is mirror-symmetric to the procedure described above, with respect to the target point $x_{i+1/2}$.

\subsection{WENO approximation for H-J equations}
\label{append2}

Given th point values $\{\phi_i\}$, we want to approximate the derivative $\partial_{\xi}\phi$ at $\xi = \xi_{i}$.
%, denoted as $\partial_{\xi}\phi^{\pm}_{i}$ and $\partial_{\xi}\phi^{+}_{i}$ are left and right-sided approximations of $\partial_{\xi}\phi$ at $\xi = \xi_{i}$. 
Introduce $\Delta^{+} \phi_{k}=\phi_{k+1}-\phi_{k}$ and $\Delta^{-} \phi_{k}=\phi_{k}-\phi_{k-1}$. Then, the fifth order approximation of $\partial_{\xi}\phi^{+}_{i}$ can be obtained as follows. 

\begin{enumerate}
	\item  There are three small stencils $S^{(r)}_{i}=\{\xi_{i-2+r},\ldots,\xi_{i+1+r}\}$, $r=0,1,2$. On each stencil, we can obtain the approximation to $\partial_\xi \phi_{i}$.
	\begin{equation*}
	\begin{aligned}	
	\partial_\xi \phi_{i}^{(0)} &= \frac{1}{\Delta\xi} \left( \frac{1}{6}\phi_{i-2}-\phi_{i-1}+\frac{1}{2}\phi_{i}+\frac{1}{3}\phi_{i+1} \right),\\	
	\partial_\xi \phi_{i}^{(1)} &= \frac{1}{\Delta\xi} \left( -\frac{1}{3}\phi_{i-1}-\frac{1}{2}\phi_{i}+\phi_{i+1}-\frac{1}{6}\phi_{i+2} \right),\\	
	\partial_\xi \phi_{i}^{(2)} &= \frac{1}{\Delta\xi} \left( -\frac{11}{6}\phi_{i}+3\phi_{i+1}-\frac{3}{2}\phi_{i+2}+\frac{1}{3}\phi_{i+3} \right), \\
	\partial_\xi \phi_{i} ^{big} &= d_0 \partial_\xi \phi_{i}^{(0)} + d_1 \partial_\xi \phi_{i}^{(1)} + d_2 \partial_\xi \phi_{i}^{(2)} ,
	\end{aligned}	
	\end{equation*}
	with linear weights $d_0=3/10$, $d_1=6/10$ and $d_2=1/10$.
	
	\item  Get the nonlinear weights 
	\begin{equation*}
	\omega_{r}=\frac{\alpha_{r}}{\alpha_{0}+\alpha_{1}+\alpha_{2}}, \quad
	\alpha_{r}=\frac{d_r}{\left(\epsilon+IS_{r}\right)^{2}}, \quad r=0,1,2.
	\end{equation*}
	Here, %$\epsilon$ is used to prevent the denominators from becoming zero, in our computation, we shall use 
	we take $\epsilon=10^{-6}$. 
	And the smoothness indicators $IS_{r}$ are defined by
	\begin{equation*}
	\begin{aligned} 
	IS_{0} &=13(b-a)^{2}+3(3b- a)^{2}, \\ 
	IS_{1} &=13(c-b)^{2}+3( c+ b)^{2}, \\ 
	IS_{2} &=13(d-c)^{2}+3( d-3c)^{2} ,
	\end{aligned}
	\end{equation*}
	and the parameters have form as
	\begin{equation*} 
	a =\frac{1}{ \Delta \xi} \Delta^+ \Delta^- \phi_{i-1}, \quad
	b =\frac{1}{ \Delta \xi} \Delta^+ \Delta^- \phi_{i  }, \quad
	c =\frac{1}{ \Delta \xi} \Delta^+ \Delta^- \phi_{i+1}, \quad
	d =\frac{1}{ \Delta \xi} \Delta^+ \Delta^- \phi_{i+2}.
	\end{equation*}
	
	\item Finally, we can obtain the WENO approximation
	\begin{equation*}
	\partial_\xi  \phi^+_{i}=\sum^2_{r=0}\omega_{r} {\partial_\xi \phi_{i}}^{(r)}.
	\end{equation*}
\end{enumerate}

The process to obtain $\partial_{\xi}\phi^{-}_{i}$ is mirror-symmetric to that for $\partial_{\xi}\phi^{+}_{i}$ with respect to the target point $\xi_{i}$.

\bibliographystyle{plain}
\bibliography{ref}	

\end{document}